\documentclass[review]{elsarticle}
\usepackage[dvipsnames]{xcolor}
\usepackage{amssymb,amsthm,amsmath,stmaryrd}
\usepackage[margin=2.5cm, bottom=2cm]{geometry}
\usepackage{subcaption}
\usepackage{algorithm, algorithmic}
\usepackage[colorlinks=true, linkcolor=blue, citecolor=red, urlcolor=blue]{hyperref}
\usepackage{todonotes}
\usepackage{lineno}
\journal{}

\newcommand{\Figref}[1]{Fig.~\ref{#1}}

\newcommand{\Sref}[1]{Section~\ref{#1}}

\newcommand{\ten}[1]{\mathcal{#1}}

\begin{document}
\begin{frontmatter}

\title{A Tensor Train-Based Isogeometric Solver for Large-Scale 3D Poisson Problems on Complex Geometries}

\author[LANL]{Quoc Thai Tran\corref{cor1}}
\author[LANL]{Duc P. Truong}
\author[LANL]{Kim \O. Rasmussen}
\author[LANL]{Boian Alexandrov}
\address[LANL]{Theoretical Division, Los Alamos National Laboratory, Los Alamos, NM, USA}
\cortext[cor1]{Corresponding author\\
Email addresses: thai.tran@lanl.gov}

\begin{abstract}
We introduce a three dimensional (3D) fully tensor train (TT)-assembled isogeometric analysis (IGA) framework, TT-IGA, for solving partial differential equations (PDEs) on complex geometries. Our method reformulates IGA discrete operators into TT format, enabling an efficient compression and computation while retaining the geometric flexibility and accuracy. Unlike prior low-rank approaches that often rely on structured domains, our framework accommodates general 3D geometries through low-rank TT representations of both the geometry mapping and the PDE's discretization. We demonstrate the effectiveness of the proposed TT-IGA framework on the three-dimensional Poisson equation, achieving substantial reductions in memory and computational cost without compromising solution quality.

\end{abstract}

\end{frontmatter}

\section{Introduction}

High-fidelity numerical simulations of high-dimensional partial differential equation (PDE) are essential across engineering and scientific applications, but achieving accuracy on complex 3D geometries often collides with computational cost. Conventional finite element methods must resolve fine-scale geometric and physical features via mesh refinement, leading to systems with billions of degrees of freedom (DOFs) and to the curse of dimensionality \cite{bellman1966dynamic}. This quickly becomes prohibitive in terms of memory and computational power, particularly for repeated simulations in optimization, control, or uncertainty quantification.

Tensor network (TN) methods, which decompose high-dimensional arrays into networks of smaller, interconnected tensors, 
have shown great promise for mitigating the curse of dimensionality and efficiently solving PDEs. These methods include formats such as tensor train (TT), canonical polyadic decomposition (CPD), and matrix product states (MPS), with the TT format~\cite{oseledets2011tensor} proving particularly effective due to its chain-like structure that enables linear scaling with dimension rather than exponential growth. Accurate and scalable TT-solvers for compressible Euler~\cite{danis_tensor-train_2025}, shallow water~\cite{danis2025SW}, Maxwell~\cite{manzini2023tensor}, Navier-Stokes~\cite{kornev2023numerical}, as well as for neutron transport equations~\cite{truong2024NTE,ortega2025NucRPRob} and incompressible fluid dynamics \cite{gourianov2022quantum}, which significantly reduce memory usage while preserving accuracy have been reported. TN methods have also been applied to kinetic models such as the Vlasov–Maxwell system~\cite{allmann2022parallel} and Vlasov–Poisson equations via quantum-inspired MPS representations~\cite{ye2022quantum}. 
Beyond direct PDE solvers, TNs have been integrated into model reduction approaches, including TT-DMD~\cite{li2023tensor} and TT-OpInf~\cite{danis_TTOpInf_2025}, as well as used for ultra fast calculations of high-dimensional integration \cite{alexandrov2023challenging}, and even calculations of the configurational integral in statistical mechanics~\cite{truong2025breaking}. These developments demonstrate TNs as a promising framework for multi-physics problems involving convection, diffusion, and transport. However, these successes have largely been confined to problems discretized on regular, Cartesian domains—a constraint that severely limits applicability to realistic engineering geometries.

While TN-based techniques have successfully captured the low-rank structure of differential operators in settings like finite differences, spectral collocation, and finite elements~\cite{adak2025tensor,manzini2023tensor}, their extension to non-Cartesian domains remains limited~\cite{markeeva2021qtt,mantzaflaris2017low,kornev2024tetrafem,ion2022tensor}. This restriction poses a major challenge for realistic simulations, where complex geometries—such as curved boundaries, embedded interfaces, or anisotropic features—are often discretized using non-Cartesian grids, isogeometric analysis, or unstructured meshes.
These geometrically flexible discretizations typically destroy the separability exploited by Cartesian TN methods, leading to rapid TT-rank growth and degraded compression. As a result, the computational benefits of TN solvers diminish, creating a critical bottleneck. Addressing this challenge requires co-designing geometry-aware discretizations and rank-preserving assembly pipelines to extend the reach of TN solvers to complex, realistic domains.
Isogeometric analysis (IGA)~\cite{hughes2005isogeometric} was introduced to unify computer-aided design (CAD) and numerical simulation by using spline basis functions—such as B-splines and NURBS—for both geometry representation and solution approximation. This allows IGA to model complex geometries without requiring mesh approximation or geometry reconstruction. In addition, the use of multivariate tensor-product bases results in structured system matrices that are well-suited for low-rank tensor decomposition methods. These features make IGA a natural candidate for integration with tensor network solvers. Recent work~\cite{guo2025tensor} has combined tensor-train representations with model reduction techniques enabling real-time and parametric simulations within the IGA framework. This highlights IGA's ability to support both geometric flexibility and algebraic structure well-suited for tensor compression.

Several recent efforts have explored the integration of tensor network (TN) techniques with isogeometric analysis (IGA) to enable low-rank solvers on structured domains. Ion et al.~\cite{ion2022tensor} developed a TT-based isogeometric solver for parameter-dependent geometries by embedding geometric variations into an extended tensor domain. While effective for parametric studies, their approach assumes structured spline-based patches and represents geometry mappings approximately, which may lead to TT-rank growth. Similarly, Mantzaflaris et al.~\cite{mantzaflaris2017low} introduced a CPD-based Galerkin-IGA framework that emphasizes Kronecker-structured assembly and low-rank quadrature, but remains confined to regular tensor-product spline domains. Markeeva et al.~\cite{markeeva2021qtt} extended this line of work to 2D problems using Quantized Tensor Train (QTT) format \cite{khoromskij2011d}, demonstrating high efficiency on uniform parametric grids. More recently, Dektor and Venturi~\cite{dektor2025coordinate} proposed a coordinate flow approach to reduce tensor ranks in general settings, offering a promising direction for improving compressibility in geometrically distorted domains. Related efforts have also explored TT acceleration for boundary integral equations on complex geometries~\cite{corona2017tensor}, where low-rank tensor approximations are used to compress kernel interactions and accelerate solution of dense systems. While effective in their respective settings, these methods are tailored to integral formulations and do not address the challenges of assembling and solving volumetric PDE systems in TT format.

Despite these advances, applying TT-based solvers to general 3D geometries remains an open challenge. Existing TT-IGA formulations typically rely on separable domains or assume low-dimensional parametric variations~\cite{mantzaflaris2017low, ion2022tensor}. In most cases, key geometric components—such as the mapping, metric tensor, and Jacobian determinant—are either approximated or assembled in full, reducing the benefits of TT compression. Addressing these limitations requires a formulation that can operate directly on complex geometries while preserving low-rank structure throughout assembly and solution.

In this work, we present a fully TT-assembled isogeometric analysis (TT-IGA) framework for solving elliptic PDEs on general 3D geometries. Unlike previous approaches~\cite{ion2022tensor, mantzaflaris2017low} that approximate the geometry mapping or rely on separable spline patches, our method retains the exact NURBS-based geometry and constructs all components—geometry, metric tensor, Jacobian determinant, stiffness matrix, and load vector—directly in the TT format. We achieve this by exploiting the tensor-product structure of B-spline bases and applying TT-cross interpolation to compress the parametric dependencies of geometric and physical quantities, all within a user-prescribed accuracy tolerance. The resulting discrete operators are fully assembled in TT format and solved using Alternating Minial Energy (AMEn) TT linear solvers. We demonstrate that our TT-IGA
enables scalable, memory-efficient simulations on complex domains with up to hundreds of millions of degrees of freedom, while preserving both numerical accuracy and geometric fidelity. To the best of our knowledge, this is the first fully TT-assembled IGA solver that handles exact geometry on unstructured spline domains without requiring geometry approximation or tensor-product simplification.

Our paper is organized as follows: In \Sref{sec:iga}, we review the B-spline/NURBS basis and IGA discretization. \Sref{sec:PoissonIGA} introduces the Poisson equation, its weak formulation and the IGA discretizations. \Sref{sec:TT formulation} presents the preliminaries on TT format and develops our TT-IGA formulation. \Sref{sec:num_examples} demonstrates numerical results. We conclude in \Sref{sec::conclude}.

\section{Isogeometric Basis Functions}
\label{sec:iga}
\subsection{Univariate B-spline Basis Functions}

Univariate B-spline basis functions are constructed as a sequence of piecewise polynomial segments defined over a non-decreasing knot vector 
\begin{equation}
\Xi = \{\xi_1, \xi_2, \ldots, \xi_{n+p+1}\}.    
\end{equation}
These basis functions are generated recursively, starting from the zeroth-order functions defined by:
\begin{equation}
N_i^0(\xi) =
\begin{cases}
1, & \text{if } \xi_i \leq \xi < \xi_{i+1}, \\
0, & \text{otherwise}.
\end{cases}    
\end{equation}
Higher-order basis functions of degree \( p \) are computed using the Cox--de Boor recursion formula:
\begin{equation}
N_i^p(\xi) = \frac{\xi - \xi_i}{\xi_{i+p} - \xi_i} N_i^{p-1}(\xi) + \frac{\xi_{i+p+1} - \xi}{\xi_{i+p+1} - \xi_{i+1}} N_{i+1}^{p-1}(\xi).    
\end{equation}
The first derivative of the B-spline basis function with respect to the parametric coordinate \( \xi \) is given by:
\begin{equation}
\frac{\partial N_i^p(\xi)}{\partial \xi} =
\frac{p}{\xi_{i+p} - \xi_i} N_i^{p-1}(\xi) - \frac{p}{\xi_{i+p+1} - \xi_{i+1}} N_{i+1}^{p-1}(\xi).    
\end{equation}
A B-spline curve defined by a set of control points \( \mathbf{P}_i \) is represented as:
\begin{equation}
\mathbf{C}(\xi) = \sum_{i=1}^m N_i^p(\xi) \mathbf{P}_i.    
\end{equation}
\subsection{Multivariate B-spline Functions and Tensor Product Constructions}

Multivariate B-spline basis functions are constructed through the tensor product of univariate B-splines. The multivariate basis is expressed as:
\begin{equation}
\bar{\mathbf{N}}^{p_1 \cdots p_N}(\xi^1, \ldots, \xi^N) = \mathbf{N}^{p_1}(\xi^1) \otimes \cdots \otimes \mathbf{N}^{p_N}(\xi^N),
\end{equation}
where the associated knot vectors span the tensor product space:
\begin{equation}
\Xi = \Xi^1 \otimes \cdots \otimes \Xi^N = \{\xi^1_1, \ldots, \xi^1_{n_1+p_1+1}\} \otimes \cdots \otimes \{\xi^N_1, \ldots, \xi^N_{n_N+p_N+1}\}.
\end{equation}

A B-spline surface can be defined as:
\begin{equation}
S(\xi_1, \xi_2) = \sum_{i=0}^n \sum_{j=0}^m N_i^p(\xi_1) N_j^q(\xi_2) \mathbf{P}_{i,j},
\end{equation}
where \( N_i^p \) and \( N_j^q \) are univariate B-spline basis functions of degree \( p \) and \( q \), respectively.

\subsection{NURBS Basis Functions and Surfaces}

\begin{figure}
\begin{center}
\minipage{0.60\textwidth}
\includegraphics[width=0.99\textwidth]{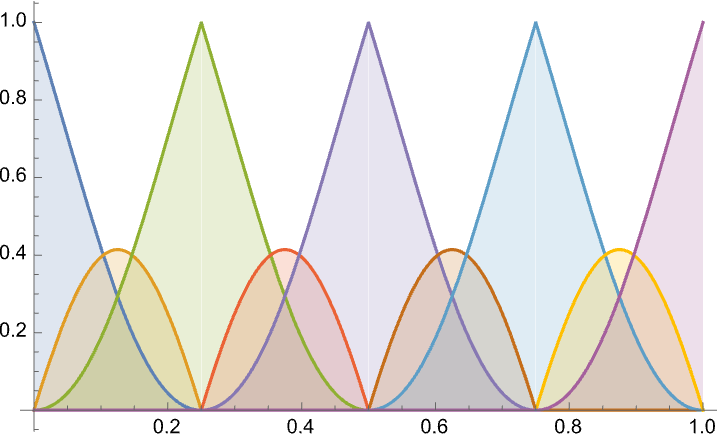}
\caption*{(a) Quadratic NURBS basis functions defined by the knot vector 
$\Xi=\{0,0,0,\,0.25,0.25,\,0.5,0.5,\,0.75,0.75,\,1,1,1\}$ and the weight vector 
$\mathbf{w}=\{1,\,1/\sqrt{2},\,1,\,1/\sqrt{2},\,1,\,1/\sqrt{2},\,1,\,1/\sqrt{2},\,1\}$.}
\endminipage
\hfill
\minipage{0.35\textwidth}
\includegraphics[width=0.99\textwidth]{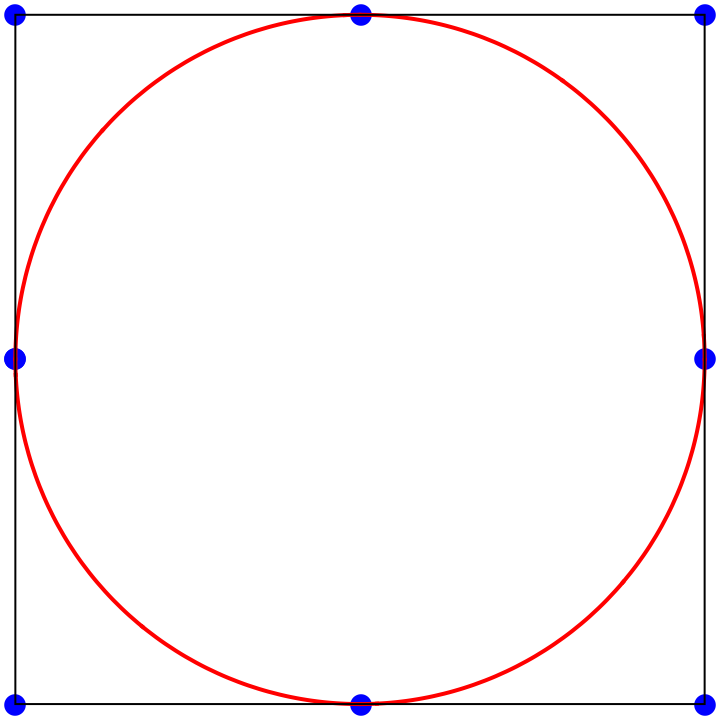}
\caption*{(b) Control points (blue dots) and the NURBS-based circular curve (red line) constructed from $\Xi$, $\mathbf{w}$, and the control points.}
\endminipage
\caption{Construction of a circle using a quadratic NURBS curve specified by the knot vector $\Xi$, the weight vector $\mathbf{w}$, and the control points.}
\label{Fig.NURBS_example}
\end{center}
\end{figure}

Non-Uniform Rational B-Splines (NURBS) extend B-spline basis functions by introducing a rational weight formulation. A univariate NURBS basis function is defined as:
\begin{equation}
R_i^p(\xi) = \frac{N_i^p(\xi) w_i}{\sum_{j=0}^n N_j^p(\xi) w_j},    
\end{equation}
where $w_i$ are the associated weights, and the knot vector $\Xi$ includes repeated knots at the ends with multiplicity $p+1$. The derivative of a univariate NURBS basis function is given by:
\begin{equation}
\frac{\partial R_i^p}{\partial \xi} =
\frac{ \frac{\partial N_i^p}{\partial \xi} W(\xi) - N_i^p(\xi) w_i \sum_{j=0}^n \frac{\partial N_j^p}{\partial \xi} w_j }{W(\xi)^2},
\quad \text{where } W(\xi) = \sum_{j=0}^n N_j^p(\xi) w_j.    
\end{equation}
A NURBS curve constructed from control points \( \mathbf{P}_i \) is written as:
\begin{equation}
\mathbf{C}(\xi) = \frac{\sum_{i=0}^n N_i^p(\xi) w_i \mathbf{P}_i}{\sum_{i=0}^n N_i^p(\xi) w_i}.    
\end{equation}
An example of constructing a NURBS based circle is presented in \Figref{Fig.NURBS_example}.

A NURBS surface is similarly defined as:
\begin{equation}
S(\xi_1, \xi_2) = \frac{\sum_{i=0}^n \sum_{j=0}^m N_i^p(\xi_1) N_j^q(\xi_2) w_{i,j} \mathbf{P}_{i,j}}{\sum_{i=0}^n \sum_{j=0}^m N_i^p(\xi_1) N_j^q(\xi_2) w_{i,j}},    
\end{equation}
with corresponding rational basis functions:
\begin{equation}
R_{i,j}^{p,q}(\xi_1, \xi_2) = \frac{N_i^p(\xi_1) N_j^q(\xi_2) w_{i,j}}{\sum_{i=0}^n \sum_{j=0}^m N_i^p(\xi_1) N_j^q(\xi_2) w_{i,j}}.    
\end{equation}
Refinement of B-spline and NURBS functions can be achieved through knot insertion, degree elevation, and continuity enhancement. Among these, knot insertion is the most commonly employed strategy. For a comprehensive treatment, refer to Piegl and Tiller~\cite{piegl1997nurbs}. In this work, we use h-refinement to refine the mesh.
\section{Poisson equation and IGA approximation}
\label{sec:PoissonIGA}

\subsection{Model Problem}
We consider the 3D Poisson equation defined over a bounded domain $\Omega \subset \mathbb{R}^3$ with appropriate boundary conditions. The governing equation is:

\begin{equation}
- \Delta u = f \quad \text{in } \Omega,
\end{equation}
where $u$ is the scalar field to be solved for, $f$ is a known source term, and $\Delta u = \nabla \cdot \nabla u$ is the Laplacian operator in physical coordinates. The boundary $\partial \Omega$ is decomposed into Dirichlet and Neumann parts:

\[
\begin{cases}
u = \bar{u} & \text{on } \Gamma_D \subset \partial \Omega \quad \text{(Dirichlet condition)} \\
\frac{\partial u}{\partial \mathbf{n}} = \bar{t} & \text{on } \Gamma_N = \partial \Omega \setminus \Gamma_D \quad \text{(Neumann condition)},
\end{cases}
\]
where $\mathbf{n}$ is the outward-pointing unit normal vector to the boundary, $\bar{u}$ is a prescribed value on the Dirichlet boundary $\Gamma_D$, and $\bar{t}$ is the prescribed Neumann traction on $\Gamma_N$.

\subsection{Weak Formulation}

To derive the weak formulation of the Poisson problem, we multiply the strong form of the equation by a test function $v \in V_0$, where $V_0$ is the space of functions that vanish on the Dirichlet boundary $\Gamma_D$, and integrate over the domain $\Omega$:

\begin{equation}
\int_\Omega (-\Delta u) v \, d\Omega = \int_\Omega f v \, d\Omega.
\end{equation}

Applying integration by parts (Green's identity) to the left-hand side and incorporating the boundary conditions, we obtain the following variational form:

\begin{equation}
\int_\Omega \nabla u \cdot \nabla v \, d\Omega = \int_\Omega f v \, d\Omega + \int_{\Gamma_N} \bar{t} v \, d\Gamma.
\end{equation}
The weak formulation of the Poisson equation is then stated as:

\begin{equation}
\text{Find } u \in V_g \text{ such that } 
\int_\Omega \nabla u \cdot \nabla v \, d\Omega = \int_\Omega f v \, d\Omega + \int_{\Gamma_N} \bar{t} v \, d\Gamma
\quad \forall v \in V_0,
\end{equation}
where the function spaces are defined as:
\begin{equation}
V_g = \{ u \in H^1(\Omega) \mid u = \bar{u} \text{ on } \Gamma_D \}, \quad 
V_0 = \{ v \in H^1(\Omega) \mid v = 0 \text{ on } \Gamma_D \}.
\end{equation}

\subsection{IGA discretizations}
In IGA, the same basis functions are used for both geometry representation and solution approximation. To exploit the advantages of the tensor-product structure, we use NURBS functions to define the geometry, and after applying h-refinement, employ B-spline functions—i.e., with all weight factors set to one—to discretize the solution. A similar approach is described in~\cite{mantzaflaris2017low}. Let $N_i(\boldsymbol{\xi})$ be the B-spline basis functions defined in the parameter space $\boldsymbol{\xi} = \{\xi_1, \xi_2, \xi_3\}$. The solution $u$ is approximated as:

\begin{equation}
u(\boldsymbol{\xi}) = \ \mathbf{N}^T(\boldsymbol{\xi}) \mathbf{u}^h 
\end{equation}
where $\mathbf{u}^h$ are the degrees of freedom. The geometry $\mathbf{x}(\boldsymbol{\xi}) = \{x(\boldsymbol{\xi}), y(\boldsymbol{\xi}), z(\boldsymbol{\xi})\}$ is computed as
\begin{equation}
    \mathbf{x}(\boldsymbol{\xi}) = \mathbf{N}^T(\boldsymbol{\xi}) \mathbf{P}   
\end{equation}
where $\mathbf{P} = \{\mathbf{P}_x, \mathbf{P}_y, \mathbf{P}_z\}$ are the control point coordinates and $\mathbf{N}(\boldsymbol{\xi})$ are the shape functions. In this work, unlike Ion et al.~\cite{ion2022tensor}, where the exact geometry is approximated by a tensor-product representation---which can result in a loss of accuracy for high-rank structures---we preserve the exact geometry without introducing any approximation. The geometry $\mathbf{x}(\boldsymbol{\xi})$ at quadrature points is computed exactly and efficiently by evaluating only the nonzero B-spline/NURBS basis functions, i.e., those for which $\xi \in [\xi_i, \xi_{i+p+1})$, together with their corresponding control points.

The discrete weak form becomes:
\begin{equation}
\mathbf{u}^h \int_{\Omega}\nabla \mathbf{N}^T \cdot \nabla \mathbf{N} \, d\Omega = \int_{\Omega} \mathbf{N}^T f \, d\Omega + \int_{\Gamma_N} \mathbf{N}^T \bar{t}\, d\Gamma.
\end{equation}
This leads to a linear system:

\begin{equation}
    \mathbf{K} \mathbf{u} = \mathbf{F}
\end{equation}

where:
\begin{equation}
\begin{split}
&\mathbf{K} = \int_{\Omega} \nabla \mathbf{N}^T \cdot \nabla \mathbf{N} \, d\Omega \quad \text{: Stiffness matrix}\\
&\mathbf{f}_v = \int_{\Omega} \textbf{N}^T f d\Omega \quad \text{: Volumetric load vector}\\
&\mathbf{f}_s = \int_{\Gamma_N} \textbf{N}^T \bar{t} \, d\Gamma \quad \text{: Surface load vector}.\\   
\end{split}
\end{equation}
The gradient operator of shape functions $\mathbf{N}$ is expressed as
\begin{equation}
    \nabla^T \mathbf{N} = \begin{bmatrix}
    \frac{\partial \textbf{N}}{x}\\
    \frac{\partial \textbf{N}}{y}\\
    \frac{\partial \textbf{N}}{z}\    
    \end{bmatrix}\
    = 
    \begin{bmatrix}
    \frac{\partial \textbf{N}}{\xi_1}\\
    \frac{\partial \textbf{N}}{\xi_2}\\
    \frac{\partial \textbf{N}}{\xi_3}  
    \end{bmatrix}^T
    \begin{bmatrix}
    \frac{\partial \xi_1}{\partial x} &\frac{\partial \xi_2}{\partial x} &\frac{\partial \xi_3}{\partial x}\\
    \frac{\partial \xi_1}{\partial y} &\frac{\partial \xi_2}{\partial y} &\frac{\partial \xi_3}{\partial y}\\
    \frac{\partial \xi_1}{\partial z} &\frac{\partial \xi_2}{\partial z} &\frac{\partial \xi_3}{\partial z}\\
    \end{bmatrix}
\end{equation}
The Jacobian matrix is defined as 
\begin{equation}
    \mathbf{J}(\boldsymbol{\xi}) = \frac{\partial \mathbf{x}}{\partial \boldsymbol{\xi}}
\end{equation}
and we have
\begin{equation}
    \nabla^T \mathbf{N} = 
    \begin{bmatrix}
    \frac{\partial \textbf{N}}{\xi_1}\\
    \frac{\partial \textbf{N}}{\xi_2}\\
    \frac{\partial \textbf{N}}{\xi_3}  
    \end{bmatrix}^T \mathbf{J}^{-1}
\end{equation}
The differential volume is computed under a change of variables as
\begin{equation}
dx dy dz = |\frac{d\mathbf{x}}{d\boldsymbol{\xi}}|d\xi_1 d\xi_2 d\xi_3
\qquad 
\text{ or }
\qquad 
d\Omega = |\frac{d\mathbf{x}}{d\boldsymbol{\xi}}|d\xi = J(\boldsymbol{\xi}) d\xi,
\end{equation}
where $J(\boldsymbol{\xi}) = |\mathbf{J}(\boldsymbol{\xi})|$ is the determinant of the Jacobian matrix $\mathbf{J}(\boldsymbol{\xi})$. The stiffness matrix and the volumetric load vector are respectively rewritten as
\begin{equation}
\label{eq:Kfinal}
    \mathbf{K} = \int_\Omega
    \begin{bmatrix}
    \frac{\partial \textbf{N}}{\xi_1}\\
    \frac{\partial \textbf{N}}{\xi_2}\\
    \frac{\partial \textbf{N}}{\xi_3}  
    \end{bmatrix}^T 
    \mathbf{R}(\boldsymbol{\xi})
    \begin{bmatrix}
    \frac{\partial \textbf{N}}{\xi_1}\\
    \frac{\partial \textbf{N}}{\xi_2}\\
    \frac{\partial \textbf{N}}{\xi_3}  
    \end{bmatrix}
d\xi
\end{equation}
and
\begin{equation}
\label{eq:ffinal} 
    \mathbf{f}_v = \int_\Omega
    f\mathbf{N}^T(\boldsymbol{\xi})
    J(\boldsymbol{\xi})
d\xi
\end{equation}
where
\begin{equation}
\label{eq.Rxi}
    \mathbf{R}(\boldsymbol{\xi}) = \mathbf{J}(\boldsymbol{\xi})^{-1}\mathbf{J}(\mathbf{\boldsymbol{\xi}})^{-T}
    J(\boldsymbol{\xi})
\end{equation}
The entries of the tensor \( \mathbf{R}(\boldsymbol{\xi}) \) from Eq.~\eqref{eq.Rxi} are denoted by \( R_{ij}(\boldsymbol{\xi}) \), where \( i, j \in \{1,2,3\} \).

Substituting these into the stiffness matrix expression \eqref{eq:Kfinal}, we obtain:
\begin{equation}
\label{eq:K_sum9_matrix_form}
\mathbf{K} = \sum_{i=1}^{3} \sum_{j=1}^{3} \mathbf{K}_{ij}, \quad \text{where} \quad \mathbf{K}_{ij} = \int_\Omega 
\left( \frac{\partial \mathbf{N}}{\partial \xi_i} \right)^T
R_{ij}(\boldsymbol{\xi})
\left( \frac{\partial \mathbf{N}}{\partial \xi_j} \right) \, d\boldsymbol{\xi}.
\end{equation}

\section{Tensor Train Formulation for Isogeometric Analysis}
\label{sec:TT formulation}
\subsection{Preliminaries}
This section introduces the core components of the tensor train (TT) methodology employed in our framework: the TT format, TT-cross interpolation, and TT linear solvers. These are well-established tools in the numerical tensor algebra literature, and we use them as existing algorithms. Our main contributions, discussed in Section~4.2, lie in the TT-compatible reformulation of IGA operators and the assembly pipeline for solving PDEs on complex 3D domains.

\subsubsection{Tensor Train Format}
\label{sec:tt_format}

The tensor train (TT) format is a hierarchical low-rank representation for high-dimensional arrays that enables efficient storage and computation. It was introduced by Oseledets~\cite{oseledets2011tensor} and has become a widely used format for representing tensors, vectors, and operators in scientific computing~\cite{khoromskij2012tensors}.

\noindent\textbf{TT format for tensors.}
Let \( \mathcal{X} \in \mathbb{R}^{n_1 \times \cdots \times n_d} \) be a $d$-dimensional tensor. The TT format expresses \( \mathcal{X} \) as a product of 3D cores:
\begin{equation}
\mathcal{X}(i_1, \dots, i_d) \approx \sum_{\alpha_1=1}^{r_1} \cdots \sum_{\alpha_{d-1}=1}^{r_{d-1}}
\ten{G}_1(1, i_1, \alpha_1) \,
\ten{G}_2(\alpha_1, i_2, \alpha_2) \cdots
\ten{G}_d(\alpha_{d-1}, i_d, 1),
\end{equation}
where \( \ten{G}_k \in \mathbb{R}^{r_{k-1} \times n_k \times r_k} \) are called TT-cores, and \( r_0 = r_d = 1 \). The integers \( r_k \) are the TT-ranks, and the overall storage is reduced from \( \mathcal{O}(n^d) \) to \( \mathcal{O}(d n r^2) \) when \( n_k = n \), \( r_k = r \).

\noindent\textbf{TT format for vectors.}
A long structured vector \( \mathbf{u} \in \mathbb{R}^{N} \), where \( N = n_1 \cdots n_d \), can be reshaped into a tensor \( \mathcal{U}(i_1, \dots, i_d) \) and approximated in TT format:
\begin{equation}
\mathbf{u}(i) \equiv \mathcal{U}(i_1, \dots, i_d) \approx \llbracket \ten{G}_1, \ldots, \ten{G}_d \rrbracket,
\end{equation}
where \( i \) is mapped from the multi-index \( (i_1, \dots, i_d) \) in lexicographic order. This representation is particularly useful for compressing and operating on high-resolution solution vectors and right-hand side vectors in PDE solvers~\cite{oseledets2012solution}.

\noindent\textbf{TT format for linear operators.}
A matrix \( \mathbf{K} \in \mathbb{R}^{N \times N} \), e.g., the stiffness matrix, can be reshaped into a $2d$-dimensional tensor:
\begin{equation}
\label{eq.Kten}
\mathcal{K}(i_1, \dots, i_d, j_1, \dots, j_d),    
\end{equation}
and approximated in TT-matrix format~\cite{truong2024NTE,dolgov2014alternating}:
\begin{equation}
\mathcal{K}(i_1, \dots, i_d, j_1, \dots, j_d) \approx 
\sum_{\alpha_1, \ldots, \alpha_{d-1}} 
\ten{M}_1(\alpha_0, i_1, j_1, \alpha_1)
\ten{M}_2(\alpha_1, i_2, j_2, \alpha_2)
\cdots
\ten{M}_d(\alpha_{d-1}, i_d, j_d, \alpha_d),
\end{equation}
where \( \ten{M}_k \in \mathbb{R}^{r_{k-1} \times n_k \times n_k \times r_k} \) are the matrix TT-cores.

\noindent\textbf{Notation: tensor product vs. Kronecker product.}  
We denote the tensor (outer) product of vectors as \( \circ \), which differs fundamentally from the Kronecker product \( \otimes \) used between matrices. The tensor product \( \circ \) of two vectors produces a higher-order tensor whose dimension is the sum of the input dimensions. For example, if \( \mathbf{a} \in \mathbb{R}^{n} \) and \( \mathbf{b} \in \mathbb{R}^{m} \), then \( \mathbf{a} \circ \mathbf{b} \in \mathbb{R}^{n \times m} \), which is a rank-one order-2 tensor (i.e., a matrix). Similarly, the outer product of three univariate function vectors results in a rank-one tensor of order 3. In contrast, the Kronecker product \( \otimes \) is defined between matrices and yields another matrix with block structure, without increasing the tensor order. In this work, we rely on tensor (outer) products to explicitly construct higher-dimensional integrand tensors, which are later compressed using TT decomposition.

\noindent\textbf{Connection between tensor product structure and rank-1 TT format.}  
A fundamental property of the TT format is that any tensor expressed as a tensor (outer) product of vector (or matrix) components can be represented as a rank-one TT format without requiring any decomposition. For example, consider a tensor \( \ten{F} \in \mathbb{R}^{n_1 \times n_2 \times n_3} \), defined by a tensor product of three vectors:  
\begin{equation}
\label{eq.Ften}
 \ten{F} = \mathbf{f}_1 \circ \mathbf{f}_2 \circ \mathbf{f}_3.   
\end{equation}
This structure naturally admits a TT representation with TT-ranks \( r_1 = r_2 = 1 \), where each TT-core corresponds exactly to one of the vectors \( \mathbf{f}_1 \), \( \mathbf{f}_2 \), or \( \mathbf{f}_3 \). This direct construction avoids any numerical TT decomposition. The same principle applies to rank-one TT-matrices. In our work, this property is used to build TT representations of IGA basis functions and operators without approximation when separability is present.

In the following, we focus on the practical construction and usage of TT decompositions in our solver. Both the interpolation of tensor-valued quantities  and the solution of linear systems are performed using the Alternating Minimal Energy (AMEn) method. AMEn provides a unified and efficient framework for TT-cross interpolation and TT-based linear solvers, which we describe in the next section.

\subsubsection{Alternating Minimal Energy for TT-Cross Interpolation and Linear Solver}
\label{sec:amen_solver}
In this work, we employ the \textit{Alternating Minimal Energy} (AMEn) method \cite{dolgov2014alternating} as the core numerical engine for both TT-cross interpolation and linear system solution. The AMEn algorithm is a powerful TT-based iterative method originally developed for solving linear systems in the TT format, and later adapted for function interpolation.

\vspace{1ex}
\noindent\textbf{TT-Cross Interpolation.}  
To construct low-rank TT representations of multidimensional tensors arising from the geometry map, Jacobian determinant, metric tensor, force vector, and stiffness matrix, we use AMEn cross interpolation. These quantities are not assembled in full; rather, they are sampled via function evaluations at carefully selected multi-index entries, guided by a recursive rank-adaptive algorithm based on the maximum volume principle~\cite{goreinov2010find}.

In our implementation, we use the \texttt{amen\_cross} routine from the \texttt{TT-Toolbox}~\cite{oseledets2019matlab}, which performs function-based TT interpolation with user-specified tolerance $\varepsilon$. TT-cross allows efficient approximation of high-dimensional tensors such as the Jacobian determinant \( J(\boldsymbol{\xi}) \), the metric tensor \( \mathbf{R}(\boldsymbol{\xi}) \), the force vector \( \mathbf{f}(\boldsymbol{\xi}) \), and the stiffness matrix \( \mathbf{K}(\boldsymbol{\xi}) \), without forming them explicitly.

\vspace{1ex}
\noindent\textbf{TT Linear Solvers.}
For solving the linear system
\[
\mathbf{K} \mathbf{u} = \mathbf{f},
\]
where both the stiffness matrix $\mathbf{K}$ and right-hand side $\mathbf{f}$ are in TT format, we use the \texttt{amen\_solve} routine. This solver applies the AMEn algorithm to iteratively minimize the residual $\|\mathbf{K} \mathbf{u} - \mathbf{f}\|_2$ by updating one TT-core of $\mathbf{u}$ at a time. An adaptive enrichment mechanism is used at each step to control TT-rank growth and ensure convergence.

The combined use of AMEn for interpolation and solving provides a unified and efficient computational framework, enabling us to handle full 3D isogeometric problems with hundreds of millions of degrees of freedom in compressed TT format.

\subsection{TT-IGA Formulation}
\label{subsection:TT-IGA formulation}

\subsubsection{TT Format of the Stiffness Matrix}
\label{sec:tt_stiffness_operator}

To express the stiffness matrix in a form suitable for tensor train (TT) decomposition, we begin with each component \( \mathbf{K}_{ij} \) as defined in Eq.~\eqref{eq:K_sum9_matrix_form}. The integrand of each term can be interpreted as a six-dimensional tensor \( \ten{K}_{ij} \), constructed from tensor-product derivatives of the basis functions and the geometry-dependent scalar coefficient tensor \( \ten{R}_{ij}(\boldsymbol{\xi}) \).

The multivariate basis function is then expressed as:
\[
\ten{N}(\boldsymbol{\xi}) = \mathbf{N}_1(\xi_1) \circ \mathbf{N}_2(\xi_2) \circ \mathbf{N}_3(\xi_3),
\]
with directional derivatives given by:
\[
\frac{\partial \ten{N}}{\partial \xi_i} = \mathbf{N}_1^{(i)}(\xi_1) \circ \mathbf{N}_2^{(i)}(\xi_2) \circ \mathbf{N}_3^{(i)}(\xi_3),
\]
where \( \mathbf{N}_d^{(i)} = \partial \mathbf{N}_d / \partial \xi_d \) if \( d = i \), and \( \mathbf{N}_d^{(i)} = \mathbf{N}_d \) otherwise.

The geometry-dependent coefficient \( R_{ij}(\boldsymbol{\xi}) \) is approximated in TT format using TT-cross interpolation with a prescribed accuracy \( \varepsilon \), such that:
\begin{equation}
\ten{R}^{TT}_{ij}(\xi_1, \xi_2, \xi_3) \approx 
\sum_{\alpha_1 = 1}^{r_1} \sum_{\alpha_2 = 1}^{r_2}
\ten{G}_1^{(ij)}(1, \xi_1, \alpha_1) \,
\ten{G}_2^{(ij)}(\alpha_1, \xi_2, \alpha_2) \,
\ten{G}_3^{(ij)}(\alpha_2, \xi_3, 1),
\end{equation}
where \( \ten{G}_d^{(ij)} \) are the TT-cores of \( R_{ij} \) in mode \( d \).

\paragraph{TT format of stiffness components.}  
Using the connection between tensor product and rank-1 TT, each \( \ten{K}^{TT}_{ij} \) is then computed as a sum over rank-1 TT terms:
\begin{equation}
\begin{split}
\ten{K}^{TT}_{ij} 
&\approx \sum_{\alpha_1 = 1}^{r_1} \sum_{\alpha_2 = 1}^{r_2}
\left( \int_{\Omega_1} \ten{G}_1^{(ij)}(1, \xi_1, \alpha_1) \, a_{ij}^{(1)}(\xi_1) \, d\xi_1 \right)
\circ
\left( \int_{\Omega_2} \ten{G}_2^{(ij)}(\alpha_1, \xi_2, \alpha_2) \, a_{ij}^{(2)}(\xi_2) \, d\xi_2 \right) \\
&\circ
\left( \int_{\Omega_3} \ten{G}_3^{(ij)}(\alpha_2, \xi_3, 1) \, a_{ij}^{(3)}(\xi_3) \, d\xi_3 \right),
\end{split}
\end{equation}
where \( a_{ij}^{(d)}(\xi_d) \) is a scalar contraction between univariate basis functions or their derivatives in direction \( d \), based on the pair \( (i, j) \).

The contraction functions \( a_{ij}^{(d)}(\xi_d) \) are defined as:
\begin{equation}
a_{ij}^{(d)}(\xi_d) =
\begin{cases}
\mathbf{d}(\xi_d) & \text{if } d = i = j, \\
\mathbf{m}(\xi_d) & \text{if } d \neq i \text{ and } d \neq j, \\
\mathbf{c}(\xi_d) & \text{if } d = i \neq j \text{ or } d = j \neq i,
\end{cases}
\end{equation}
with
\begin{equation}
\begin{split}
\mathbf{d}(\xi_d) &= \mathbf{N}_d'(\xi_d)^T \mathbf{N}_d'(\xi_d), \\
\mathbf{m}(\xi_d) &= \mathbf{N}_d(\xi_d)^T \mathbf{N}_d(\xi_d), \\
\mathbf{c}(\xi_d) &= \mathbf{N}_d'(\xi_d)^T \mathbf{N}_d(\xi_d).
\end{split}
\end{equation}

\paragraph{Global assembly.}  
The TT format of the stiffness operator is formed by summing all TT components:
\begin{equation}
\ten{K}^{TT} = \sum_{i=1}^{3} \sum_{j=1}^{3} \ten{K}^{TT}_{ij},
\end{equation}

As the sum over \( \ten{K}^{TT}_{ij} \) may increase the TT-rank, we apply TT-rounding (also known as TT-recompression) after each summation step to maintain a low-rank representation~\cite{oseledets2011tensor}. TT-rounding ensures that the resulting TT-cores remain compressed within a user-specified tolerance \( \varepsilon \), while avoiding unnecessary rank inflation.

Next, we will describe how to construct the TT format of the force vector.

\subsubsection{Tensor Train Representation of the Force Vector}

We aim to compute the volumetric load vector from Eq.~\eqref{eq:ffinal}:
\[
\mathbf{f}_v = \int_\Omega f(\boldsymbol{\xi}) \, \mathbf{N}^T(\boldsymbol{\xi}) \, J(\boldsymbol{\xi}) \, d\boldsymbol{\xi},
\]
where \( \mathbf{N}(\boldsymbol{\xi}) = \mathbf{N}_1(\xi_1) \circ \mathbf{N}_2(\xi_2) \circ \mathbf{N}_3(\xi_3) \) is the tensor-product basis, and \( J(\boldsymbol{\xi}) \) is the determinant of the Jacobian matrix.

To exploit the low-rank structure of the integrand, we first approximate the tensor \( \ten{FJ}=f(\boldsymbol{\xi}) J(\boldsymbol{\xi}) \) in TT format using TT-cross interpolation :
\begin{equation}
\ten{FJ}^{TT}(\xi_1, \xi_2, \xi_3) \approx 
\sum_{\alpha_1 = 1}^{r_1} \sum_{\alpha_2 = 1}^{r_2}
\ten{G}_1(1, \xi_1, \alpha_1) \,
\ten{G}_2(\alpha_1, \xi_2, \alpha_2) \,
\ten{G}_3(\alpha_2, \xi_3, 1),
\end{equation}
where \( \ten{G}_d \) are the TT-cores of \( \ten{J}^{TT} \) in dimension \( d \).

Substituting this decomposition into the integrand, and leveraging the tensor product structure of the basis functions, the TT format of the volumetric force tensor, $\ten{F}^{TT}_v$, can be computed as a sum of rank-1 TT terms:
\begin{equation}
\begin{split}
\ten{F}^{TT}_v
&\approx \sum_{\alpha_1 = 1}^{r_1} \sum_{\alpha_2 = 1}^{r_2}
\left( \int_{\Omega_1}  \ten{G}_1(1, \xi_1, \alpha_1) \, \mathbf{N}_1(\xi_1) \, d\xi_1 \right)
\circ
\left( \int_{\Omega_2} \ten{G}_2(\alpha_1, \xi_2, \alpha_2) \, \mathbf{N}_2(\xi_2) \, d\xi_2 \right) \\
&\circ
\left( \int_{\Omega_3} \ten{G}_3(\alpha_2, \xi_3, 1) \, \mathbf{N}_3(\xi_3) \, d\xi_3 \right).
\end{split}
\end{equation}

\subsection{Handling Boundary Conditions}
The Boundary condition is handled in TT-format by subtracting the boundary term from the RHS. More details about this procedure is described in our previous space-time works~\cite{adak2025tensor,dibyendu2025space}. This procedure results into a linear system in TT format of only the unknowns.

\section{Numerical Experiments}
\label{sec:num_examples}
In the first two numerical examples, we investigate the convergence behavior and quantify the numerical error of the TT-IGA method using problems with known analytical solutions to the Poisson equation. The accuracy of the numerical solution is measured by the $L_2$ norm of the error, defined as
\begin{equation}
\label{eq:L2norm}
\|e\|_{L_2} = \frac{
  \displaystyle \int_{\Omega}
    \left\| \mathbf{u} - \hat{\mathbf{u}} \right\| \, d\Omega
}{
  \displaystyle \int_{\Omega}
    \left\| \hat{\mathbf{u}} \right\| \, d\Omega
},
\end{equation}
where $\hat{\mathbf{u}}$ denotes the analytical solution, and $\mathbf{u}$ is the numerical solution. In the final numerical example, we directly compare the performance of the full-grid IGA using a sparse solver with our tensor train method, to rigorously evaluate the efficiency gains achieved by the proposed approach.
\subsection{Experiment 1: L-shape domain}
\begin{figure}
\begin{center}
\minipage{0.45\textwidth}
\includegraphics[width=0.99\textwidth]{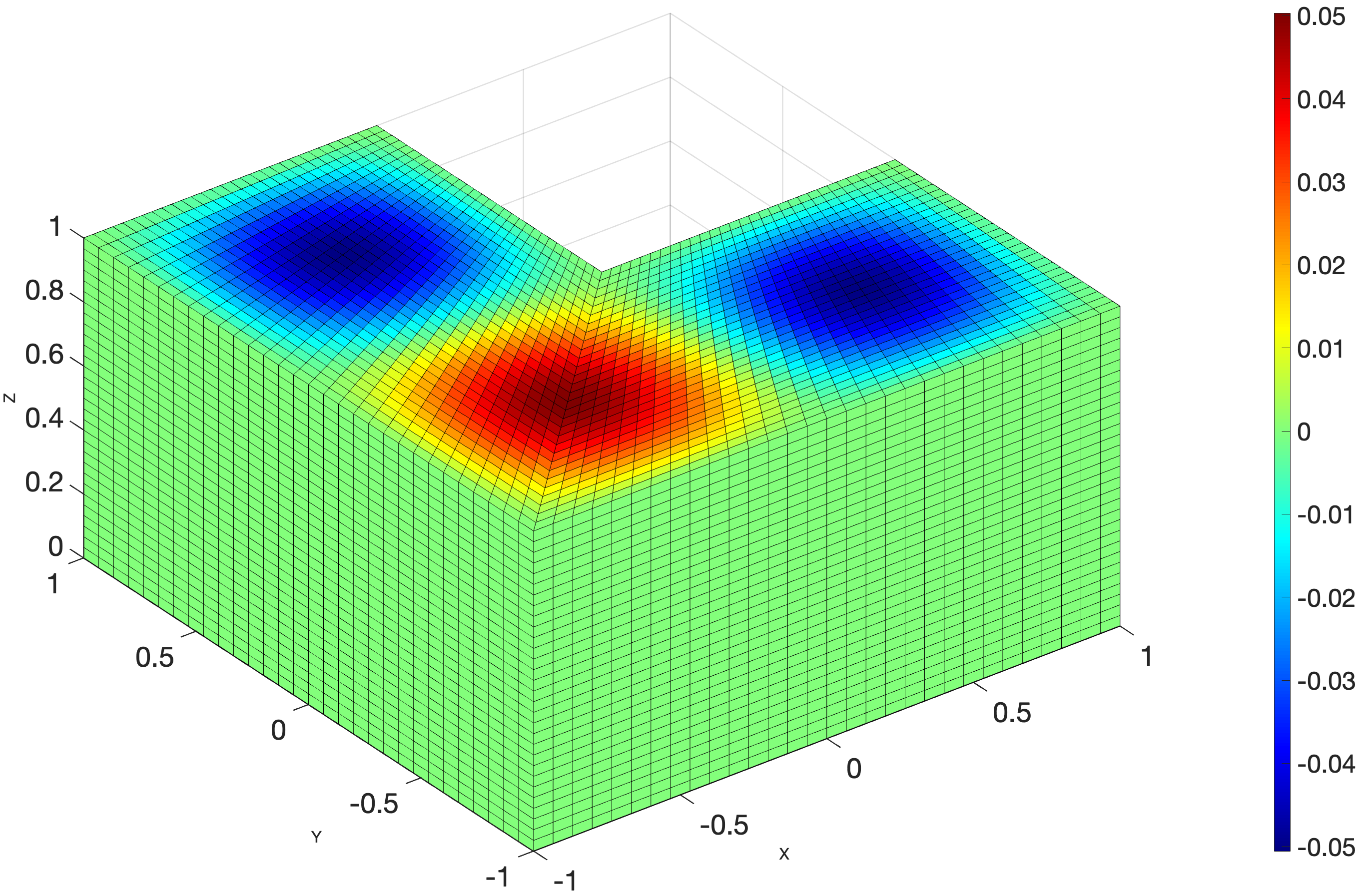}
\caption*{(a) Numerical solution}
\endminipage
\hfill
\minipage{0.45\textwidth}
\includegraphics[width=0.99\textwidth]{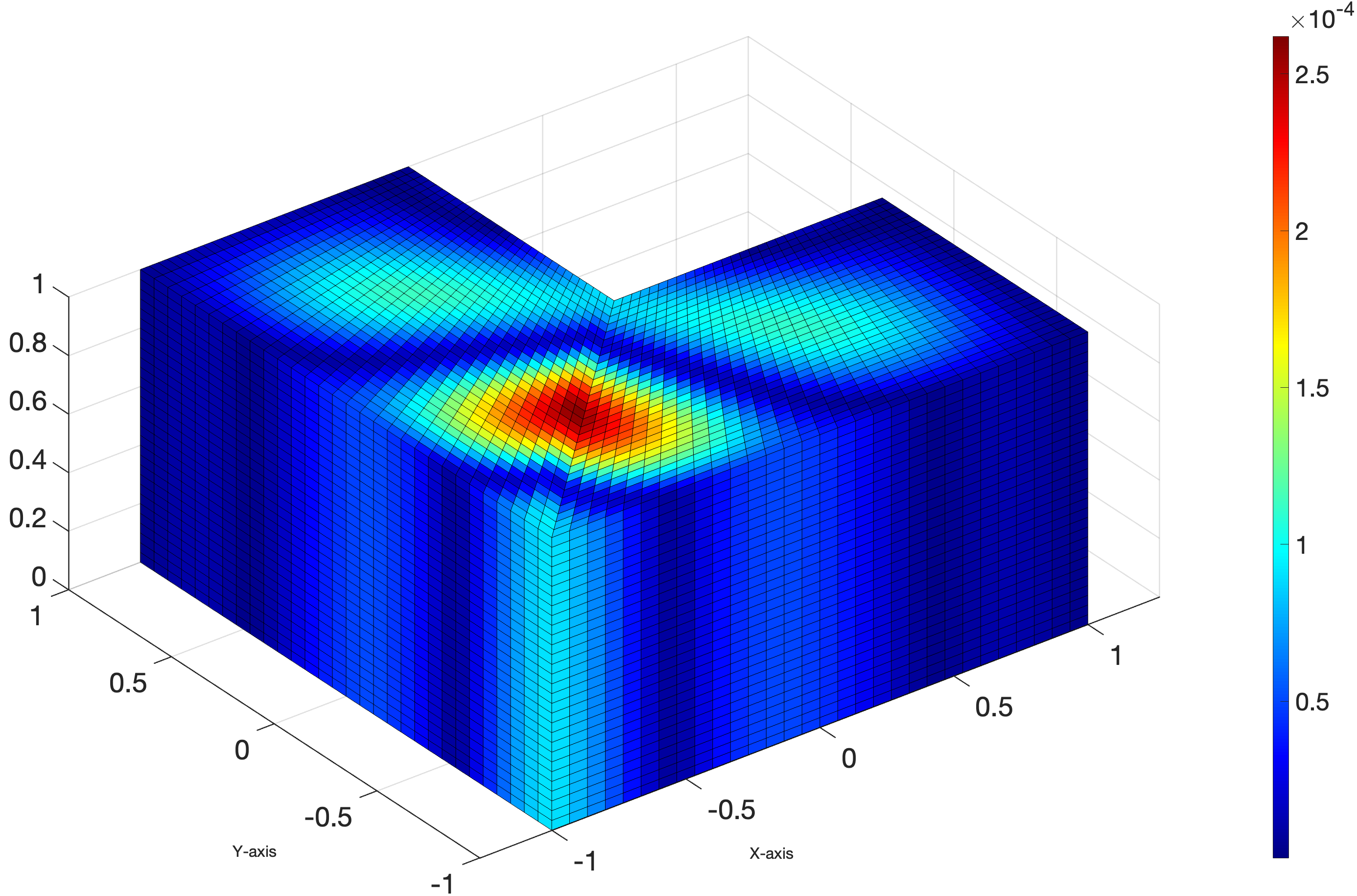}
\caption*{(b) Numerical error}
\endminipage
\caption{Numerical solution and error for the L-shaped domain with source term $f(x, y) = \sin(\pi x) \sin(\pi y)$.}
\label{Fig.Lshape_Solution_Error}
\end{center}
\end{figure}

\begin{figure}
\begin{center}
\minipage{0.5\textwidth}
\includegraphics[width=0.99\textwidth]{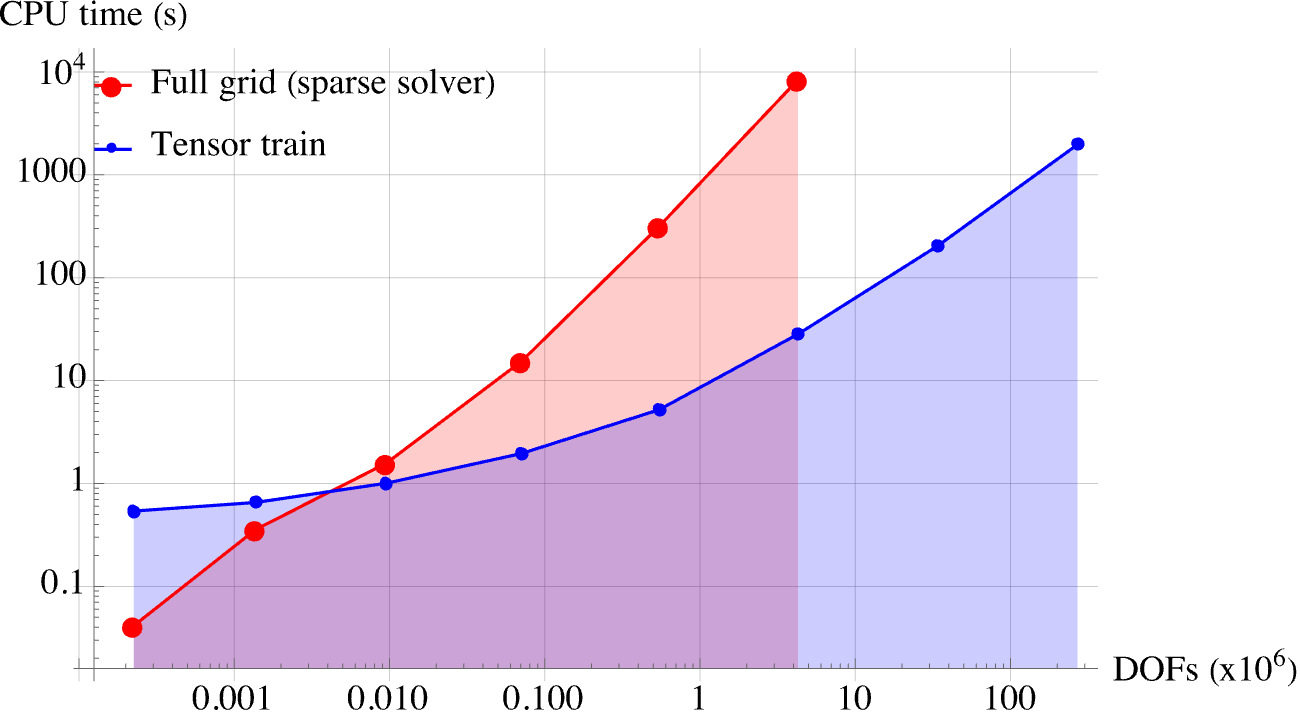}
\caption*{(a) CPU time}
\endminipage
\hfill
\minipage{0.5\textwidth}
\includegraphics[width=0.99\textwidth]{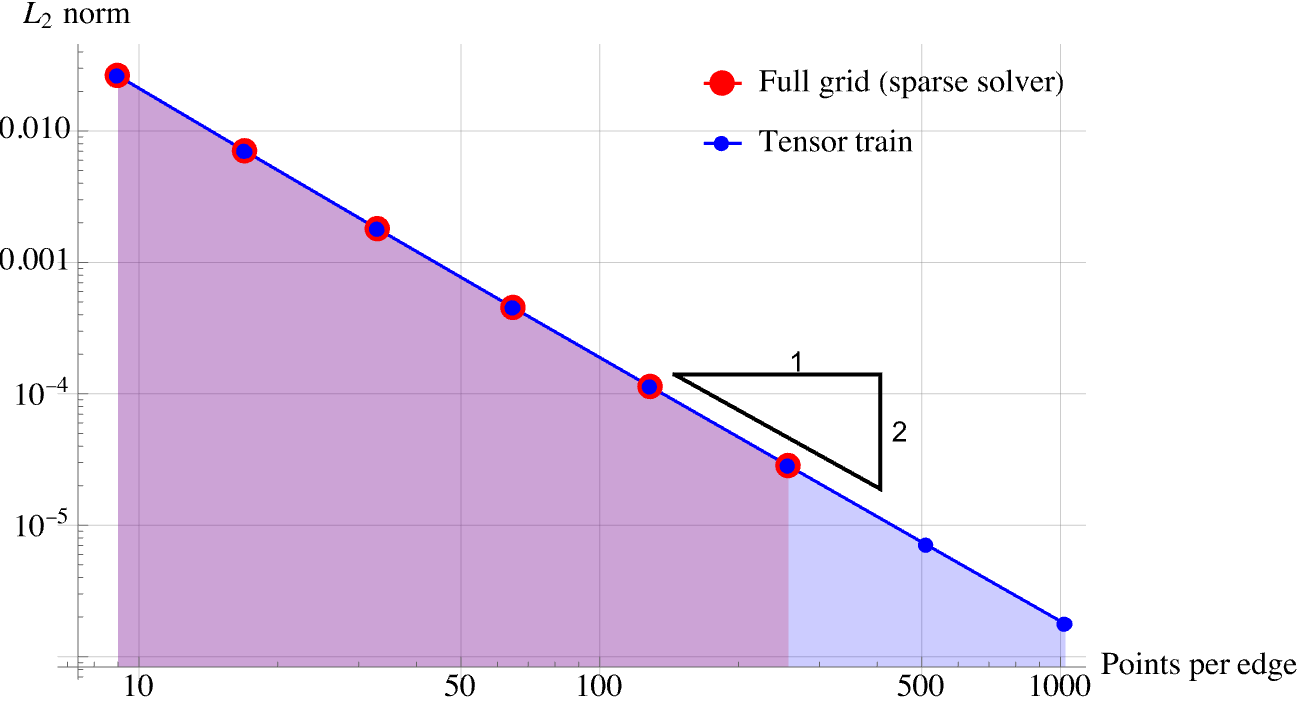}
\caption*{(b) $L_2$ norm}
\endminipage
\vfill
\minipage{0.5\textwidth}
\includegraphics[width=0.99\textwidth]{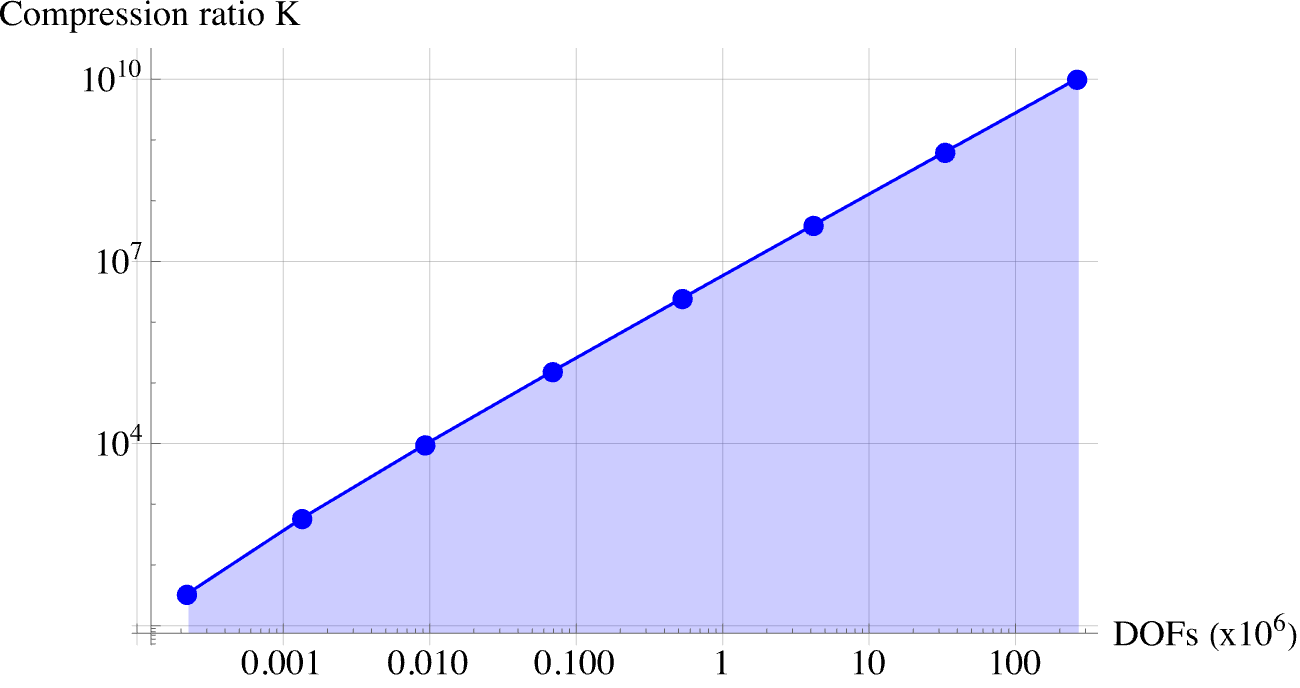}
\caption*{(c) Compression ratio of stiffness matrix \textbf{K}}
\endminipage
\hfill
\minipage{0.5\textwidth}
\includegraphics[width=0.99\textwidth]{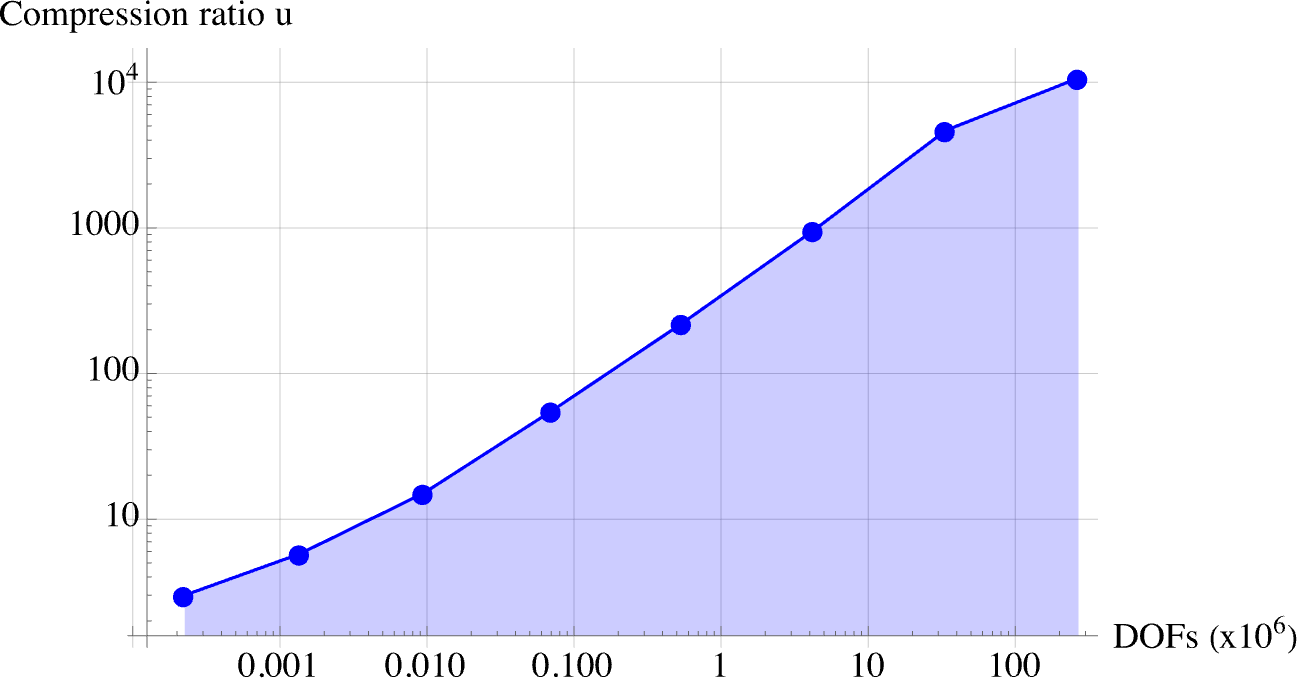}
\caption*{(d) Compression ratio of solution \textbf{u}}
\endminipage
\caption{L-shape domain - CPU time and compression ratio of stiffness matrix \textbf{K} and numerical solution \textbf{u} versus total degrees of freedom (DOFs), and $\mathrm{L_2}$ norm of the numerical solution against number of points per edge.}
\label{Fig.Lshape_convegence}
\end{center}
\end{figure}

In this example, an L-shaped domain subjected to a body force $f(x, y) = \sin(\pi x) \sin(\pi y)$ is investigated with the boundary condition $u|_{\partial \Omega} = 0$. 
The L-shaped domain is constructed by removing a smaller cube from one corner of a larger cube. The outer cube occupies the region $[-1, 1] \times [-1, 1] \times [0, 1]$, while the inner cutout has dimensions $[0, 1] \times [0, 1] \times [0, 1]$. The mesh is discretized using linear shape functions in all three parametric directions. The analytical solution to this problem is given by:
\begin{equation}
u(x, y, z) = \frac{1}{2 \pi^2} \sin(\pi x) \sin(\pi y)
\end{equation}
\Figref{Fig.Lshape_Solution_Error} illustrates the computed numerical solution and the corresponding error distribution for a mesh of size $61 \times 31 \times 31$ mesh. In \Figref{Fig.Lshape_convegence}, we report the CPU time, the $L_2$ error norm and the compression ratios for both the tangent matrix $\mathbf{K}$ and the numerical solution $\mathbf{u}$, as functions of the total degrees of freedom (DOFs).

The results demonstrate that the TT-IGA approach achieves a convergence rate consistent with the expected behavior for linear basis functions, with the slope of the $L_2$  error versus the number of points per edge approximately equal to 2. Furthermore, as the problem size increases, the TT-based solver exhibits substantially lower computational cost compared to the full-grid (sparse) approach, while maintaining excellent accuracy. The compression ratios for both the matrix and solution vector also increase with mesh refinement, confirming the efficiency and scalability of the TT-IGA method.

\subsection{Experiment 2: 3D Ring}
A 3D ring structure with inner radius $r_{\text{in}} = 0.5$, outer radius $r_{\text{out}} = 1.0$, and height $h = 1$ is studied in this example. 
Dirichlet boundary conditions are applied on the inner and outer cylindrical surfaces, with $u_{\text{in}} = 1$ at $r = r_{\text{in}}$ and $u_{\text{out}} = 2$ at $r = r_{\text{out}}$. 
The analytical solution to the Poisson equation for this problem is given by
\begin{equation}
u(r) = \frac{u_{\text{in}} \log\left(\frac{r_{\text{out}}}{r}\right) + u_{\text{out}} \log\left(\frac{r}{r_{\text{in}}}\right)}{\log\left(\frac{r_{\text{out}}}{r_{\text{in}}}\right)}.
\end{equation}
This is a homogeneous Laplace problem with \( f = 0 \), where the solution depends only on the radial coordinate \( r \) in cylindrical coordinates.

\begin{figure}
\begin{center}
\minipage{0.45\textwidth}
\includegraphics[width=0.99\textwidth]{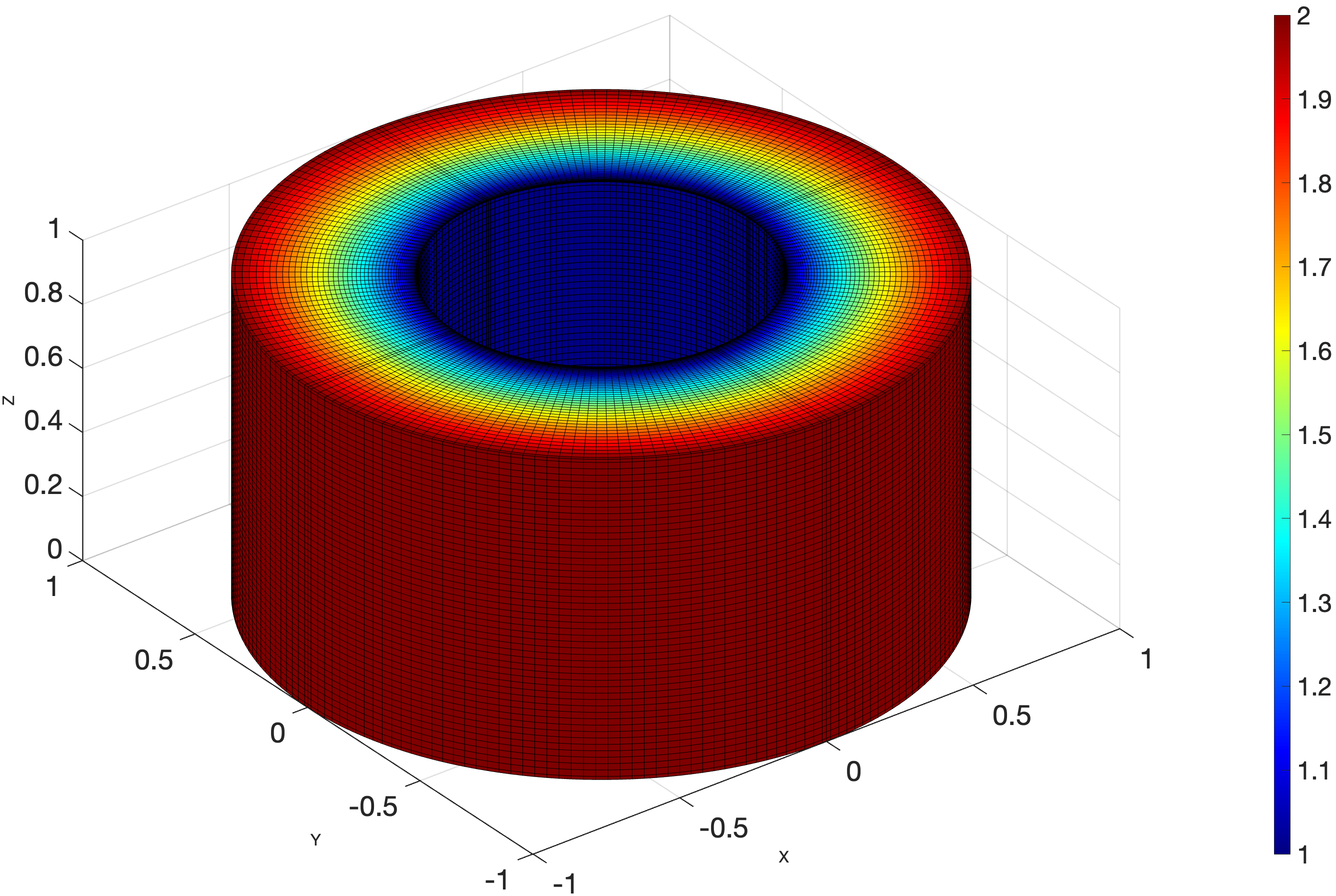}
\caption*{(a) Numerical solution}
\endminipage
\hfill
\minipage{0.45\textwidth}
\includegraphics[width=0.99\textwidth]{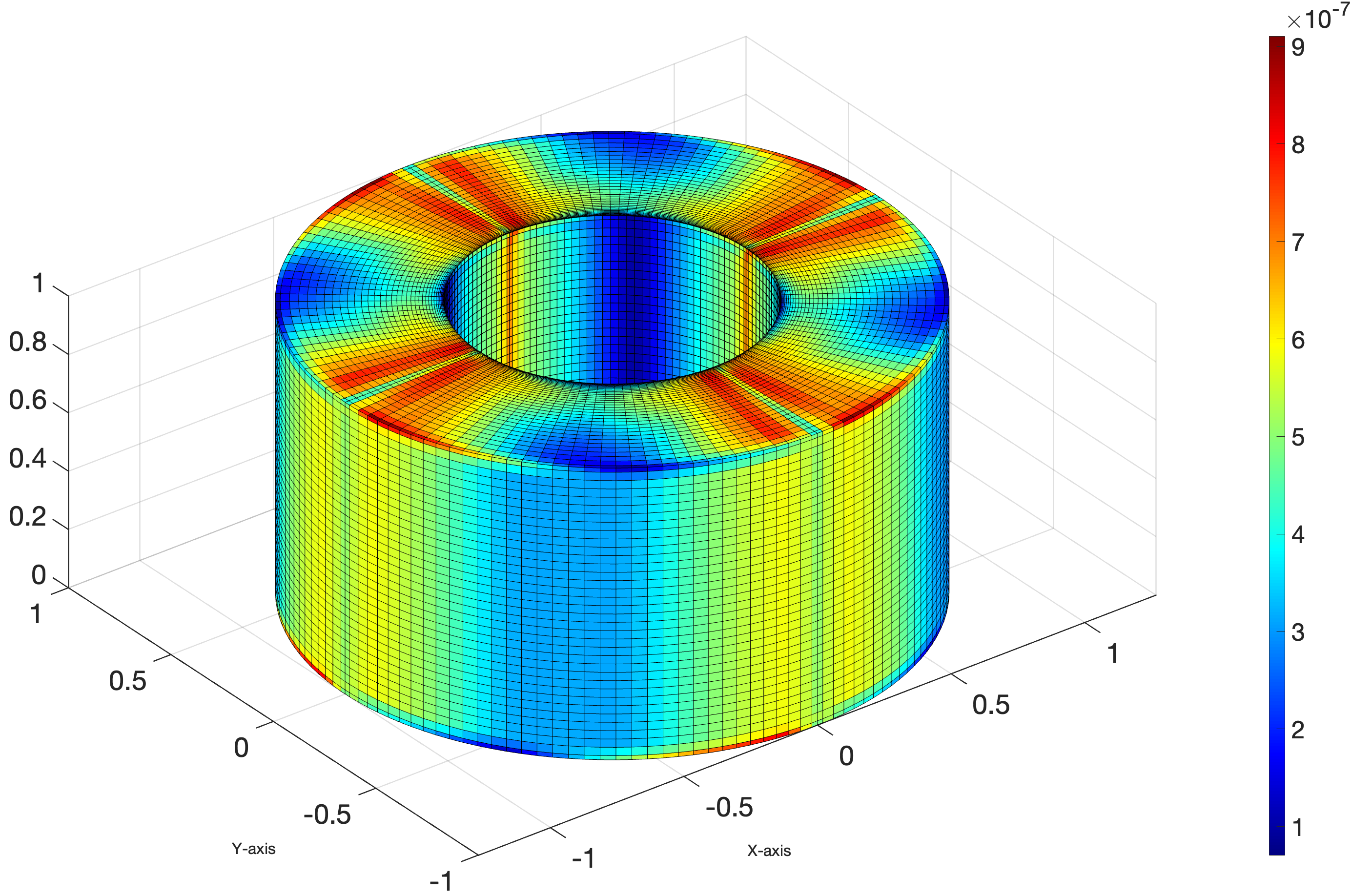}
\caption*{(b) Numerical error}
\endminipage
\caption{Numerical solution and error for the 3D ring domain with Dirichlet boundary conditions $u_{\text{in}} = 1$ and $u_{\text{out}} = 2$.}
\label{Fig.Ring_Solution_Error}
\end{center}
\end{figure}

\begin{figure}
\begin{center}
\minipage{0.5\textwidth}
\includegraphics[width=0.99\textwidth]{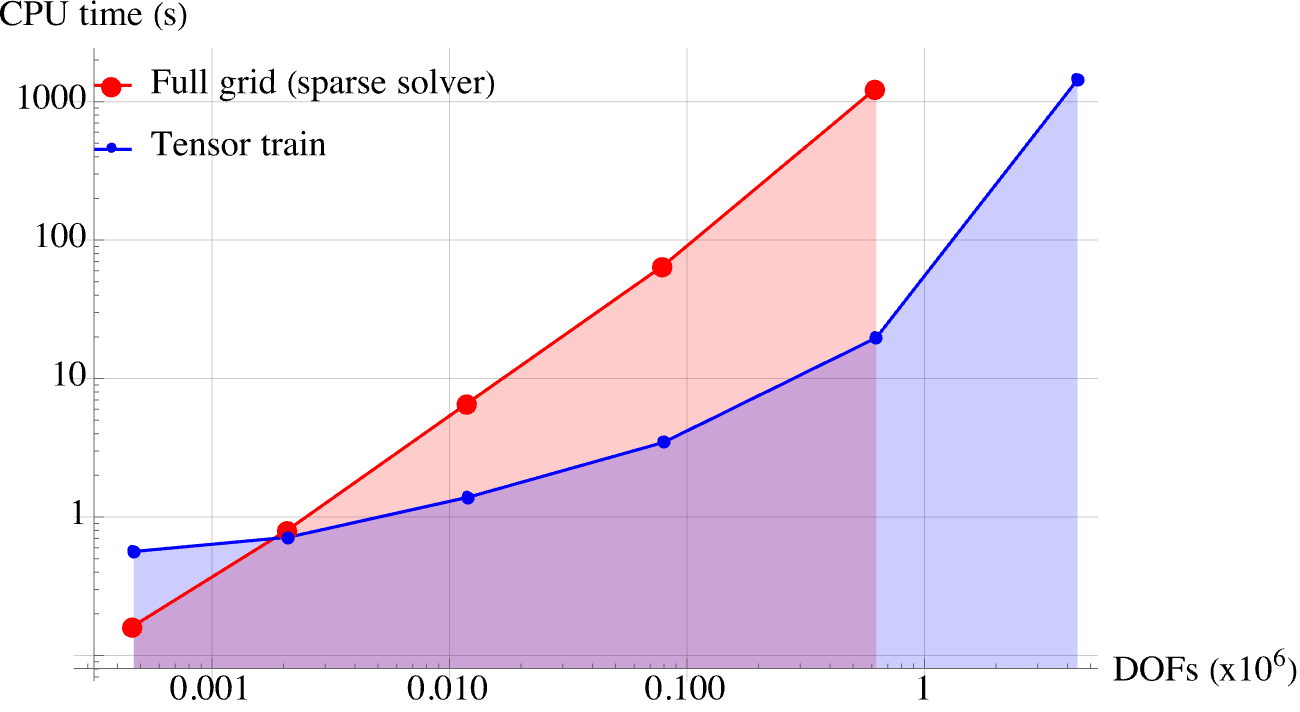}
\caption*{(a) CPU time}
\endminipage
\hfill
\minipage{0.5\textwidth}
\includegraphics[width=0.99\textwidth]{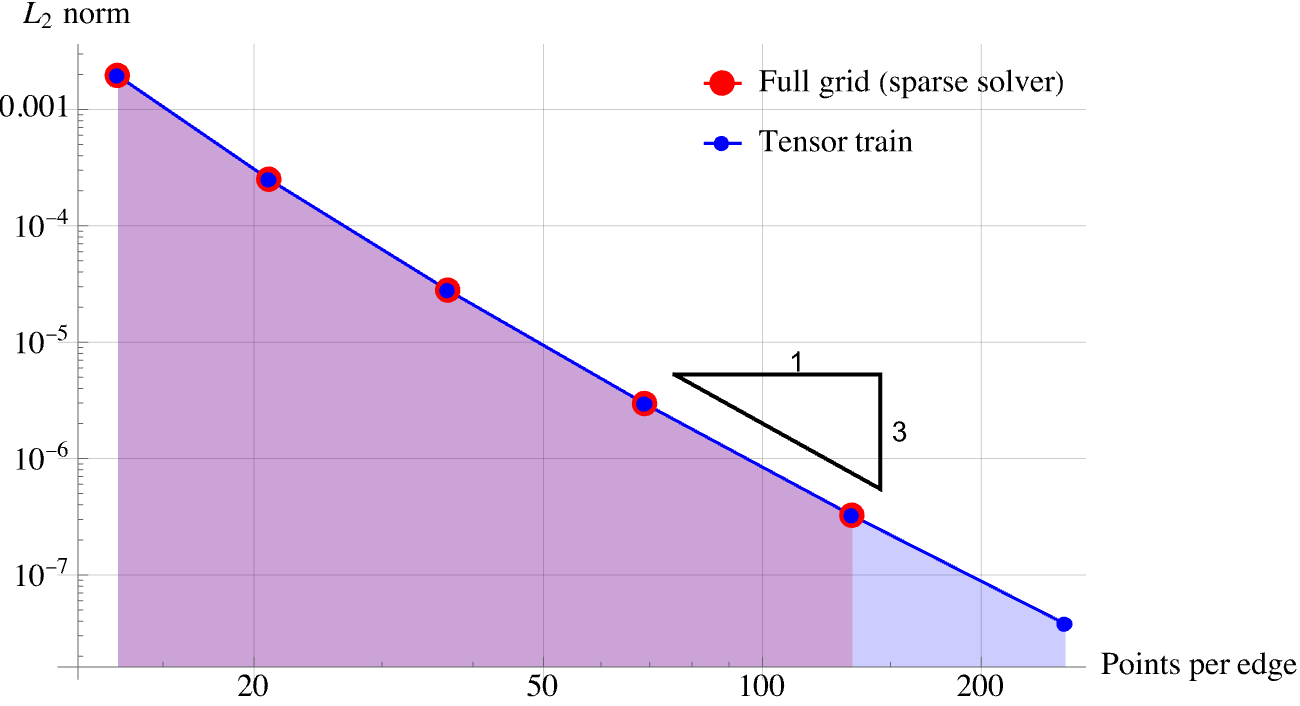}
\caption*{(b) $L_2$ norm}
\endminipage
\vfill
\minipage{0.5\textwidth}
\includegraphics[width=0.99\textwidth]{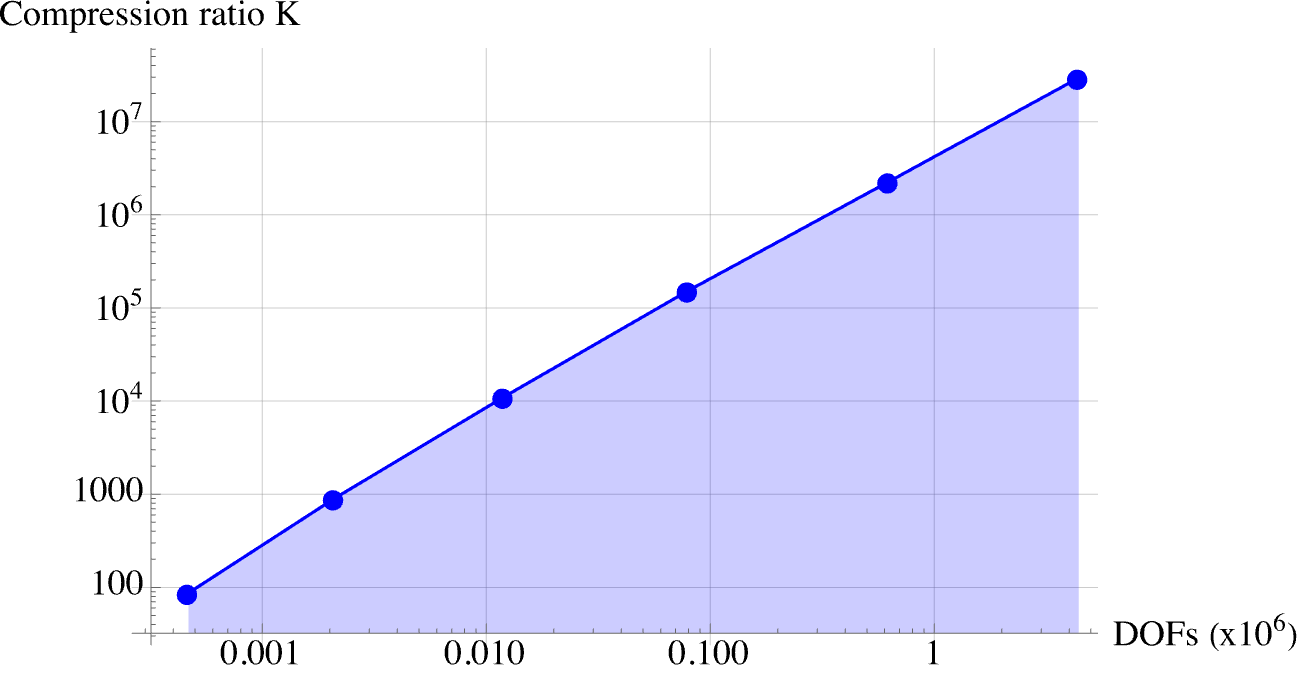}
\caption*{(c) Compression ratio of stiffness matrix \textbf{K}}
\endminipage
\hfill
\minipage{0.5\textwidth}
\includegraphics[width=0.99\textwidth]{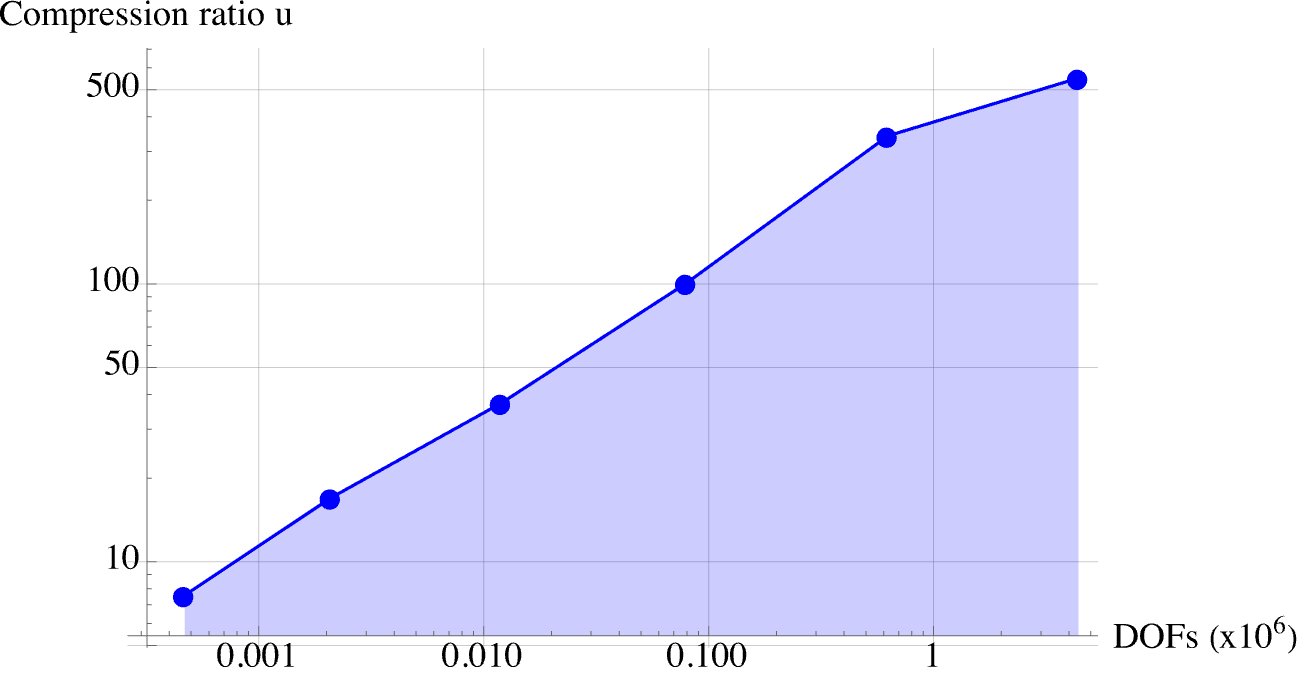}
\caption*{(d) Compression ratio of solution \textbf{u}}
\endminipage
\caption{Ring domain - CPU time and compression ratio of stiffness matrix \textbf{K} and numerical solution \textbf{u} versus total degrees of freedom (DOFs), and $\mathrm{L_2}$ norm of the numerical solution against number of points per edge.}
\label{Fig.Ring_convegence}
\end{center}
\end{figure}

\Figref{Fig.Ring_Solution_Error} presents the computed numerical solution and error distribution for a mesh size of $147 \times 37 \times 37$. The computational efficiency of the TT-IGA solver is highlighted in \Figref{Fig.Ring_convegence}, which compares CPU time and reports compression ratios for both the stiffness matrix and solution vector as the number of DOFs increases.

The convergence rate, observed as the slope of the $L_2$ error norm versus mesh resolution is approximately 3 matching the theoretical expectation for quadratic basis functions used in this experiment. As with the previous example, the TT-IGA method delivers a substantial reduction in computational time and memory requirements compared to the full-grid sparse solver.
\subsection{Experiment 3: Extended benchmarking on multiple geometries}
\begin{figure}[H]
  \centering
  \begin{subfigure}{0.3\textwidth}
    \includegraphics[width=\linewidth]{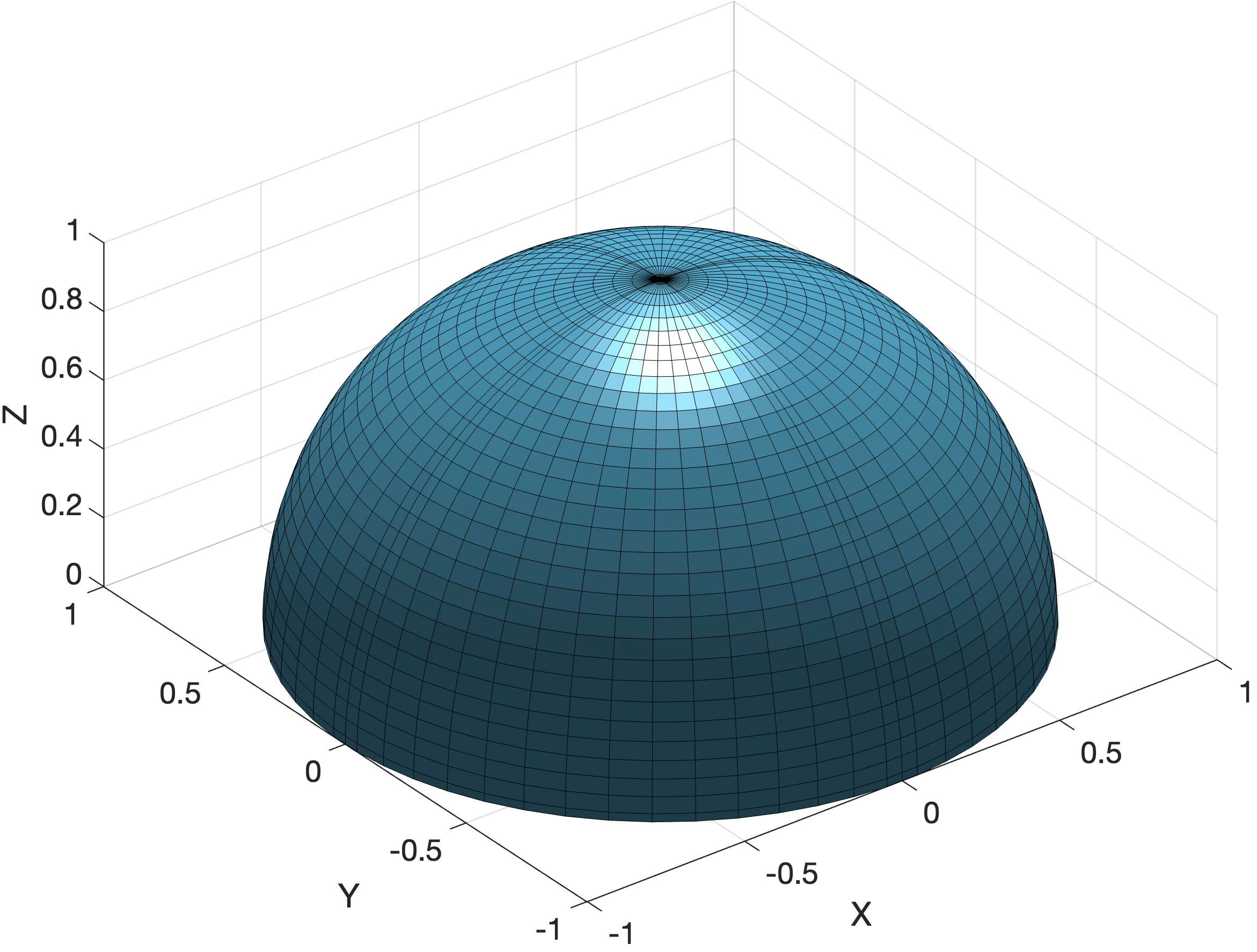}
    \caption{Closed hemisphere: $r_\text{in} = 0.5$, $r_\text{out} = 1.0$}
  \end{subfigure}
  \hfill
  \begin{subfigure}{0.3\textwidth}
    \includegraphics[width=\linewidth]{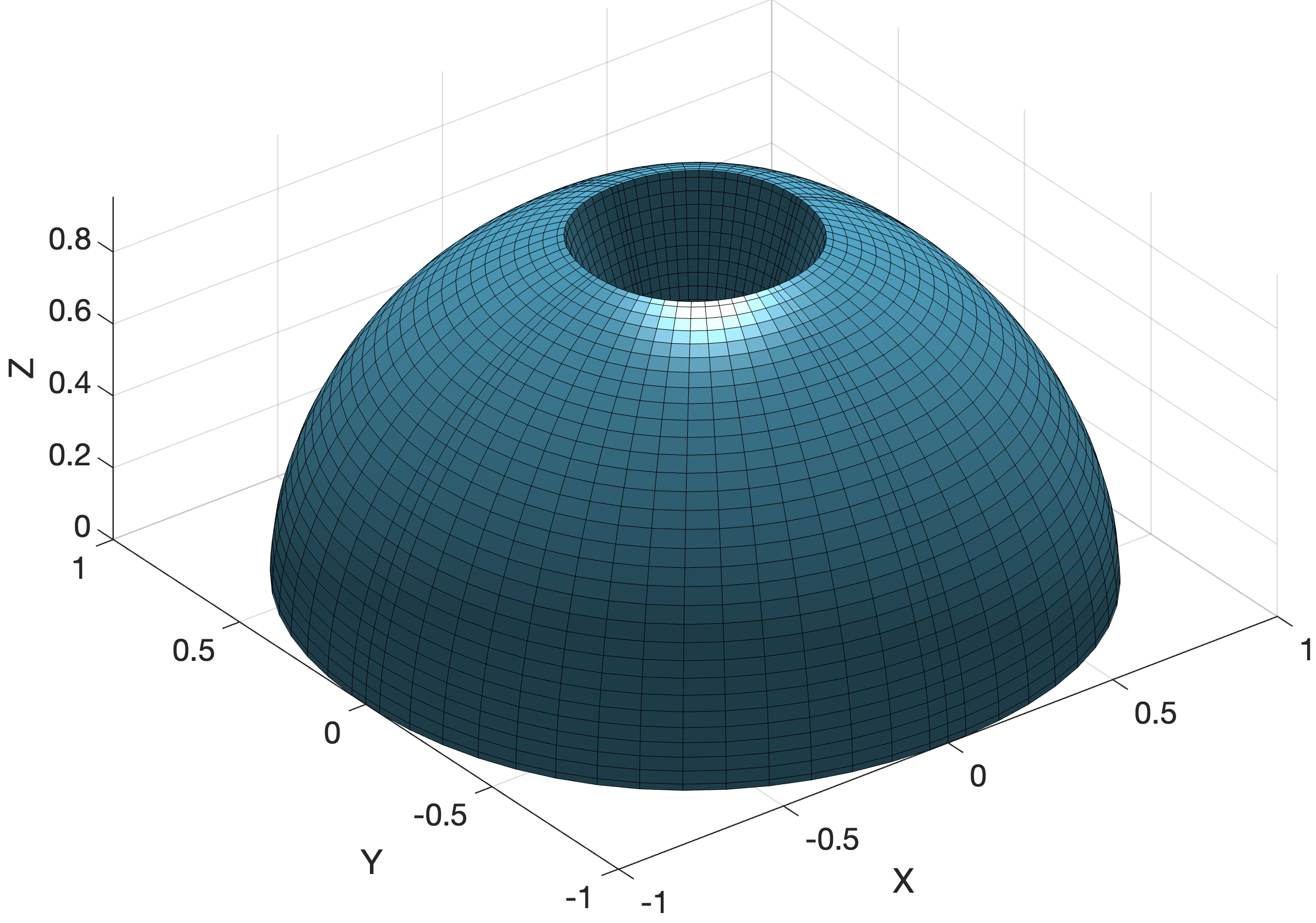}
    \caption{Opened hemisphere: $r_\text{in} = 0.5$, $r_\text{out} = 1.0$ with $18^\circ$ hole on top.}
  \end{subfigure}
  \hfill
  \begin{subfigure}{0.3\textwidth}\includegraphics[width=\linewidth]{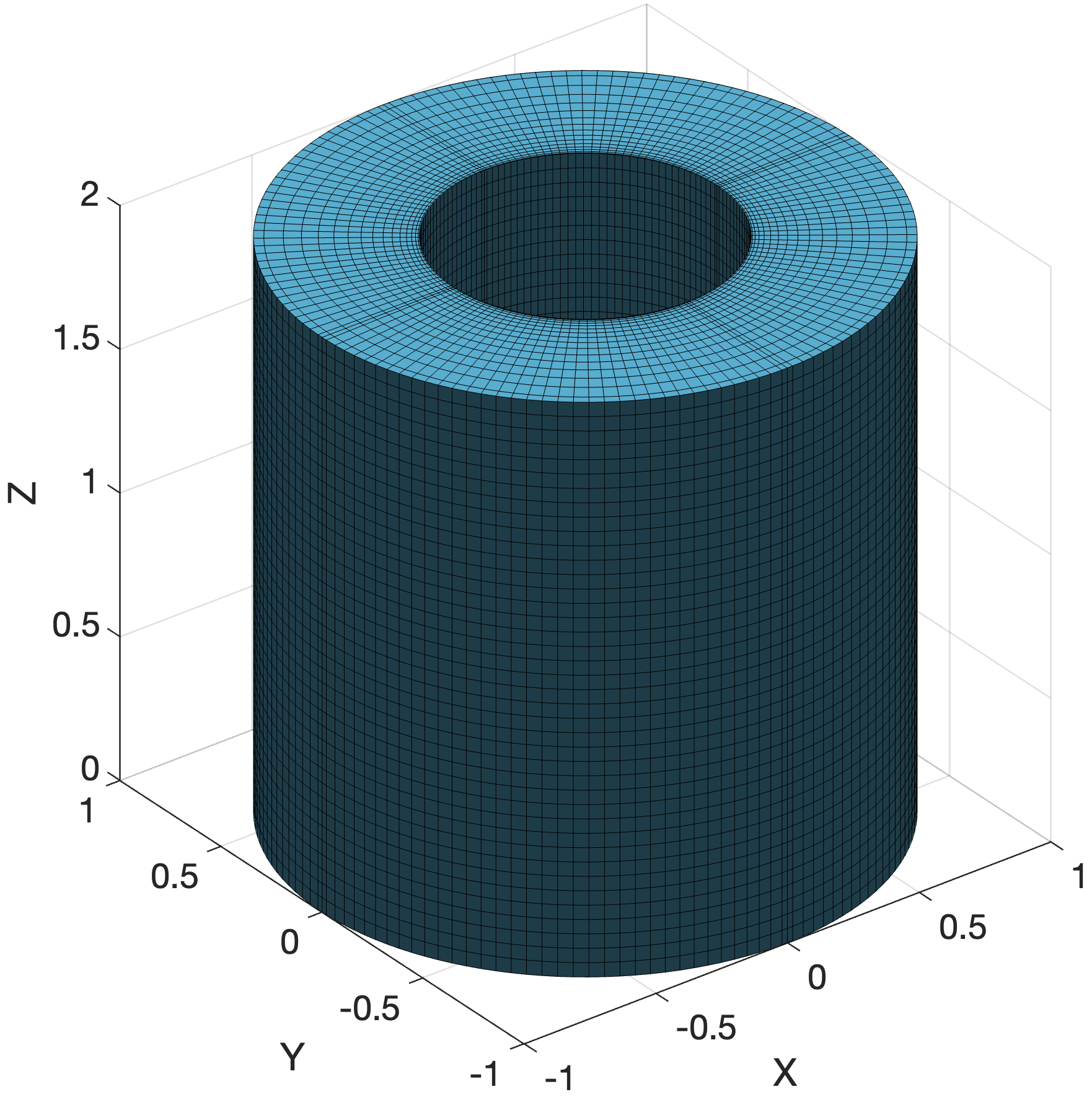}
    \caption{3D Ring: $r_\text{in} = 0.5$, $r_\text{out} = 1.0$ and $h = 2$.}
  \end{subfigure}
  \vfill
  \begin{subfigure}{0.3\textwidth}
    \includegraphics[width=\linewidth]{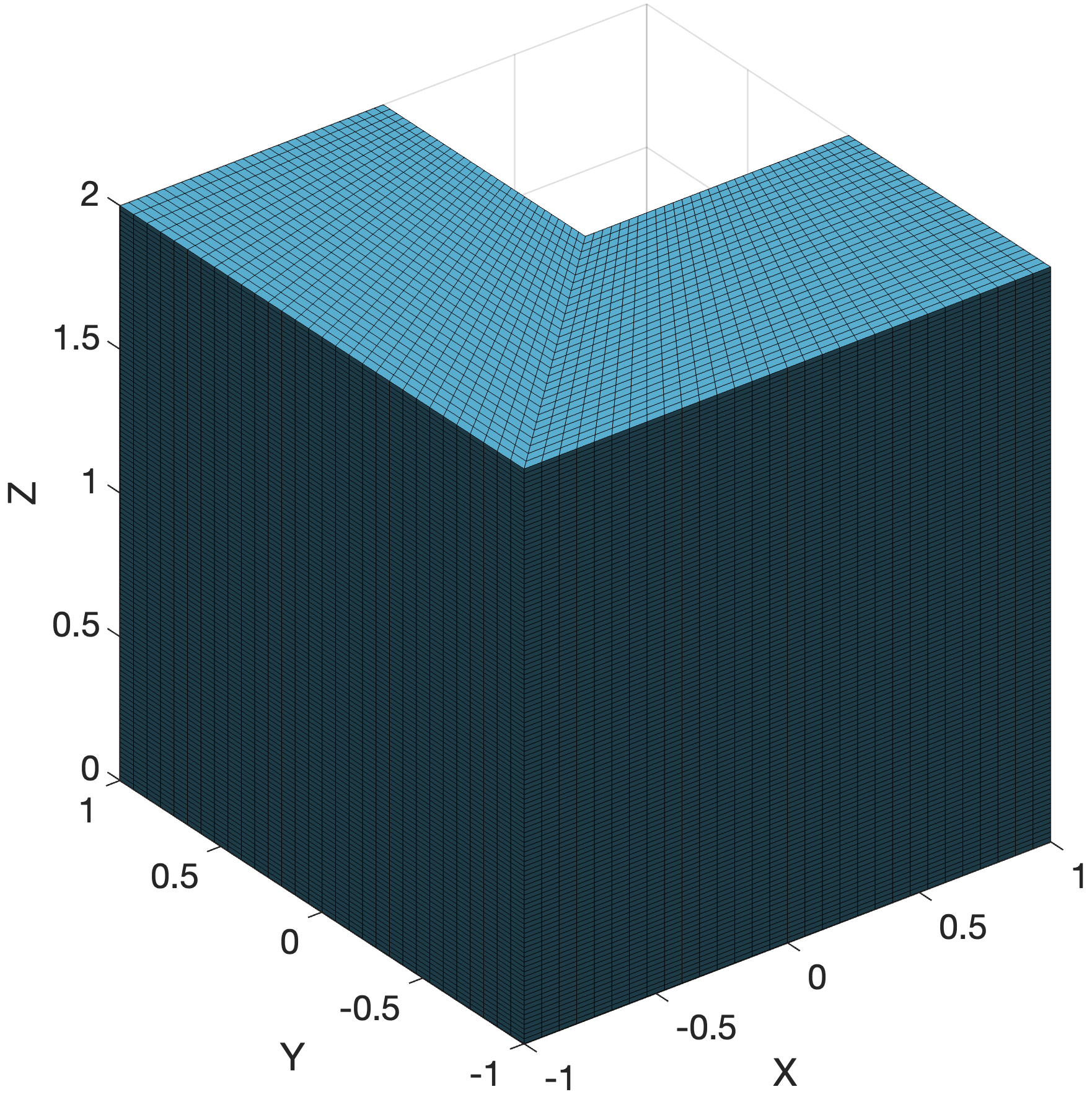}
    \caption{3D L-shape: outer cube $[-1, 1] \times [-1, 1] \times [0, 2]$ and the inner cutout $[0, 1] \times [0, 1] \times [0, 2]$.}
  \end{subfigure}
  \hfill
  \begin{subfigure}{0.3\textwidth}
    \includegraphics[width=\linewidth]{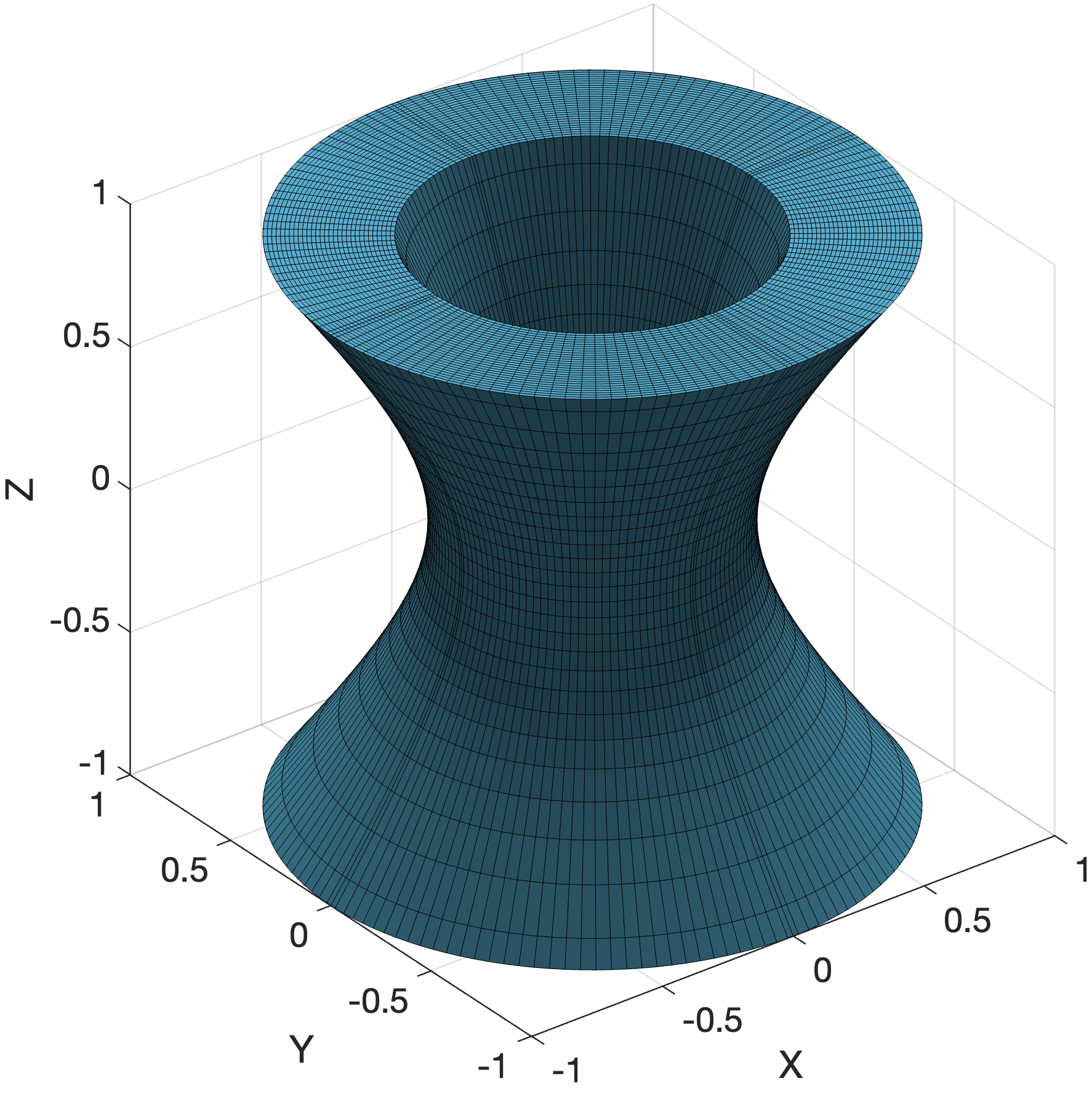}
    \caption{Hyperboloid: $r_\text{middle} = 0.5$, $r_\text{top} = 1.0$, radial thickness $r_0 = 0.3$ and $h = [-1, 1]$.}
  \end{subfigure}
  \hfill
  \begin{subfigure}{0.3\textwidth}
    \includegraphics[width=\linewidth]{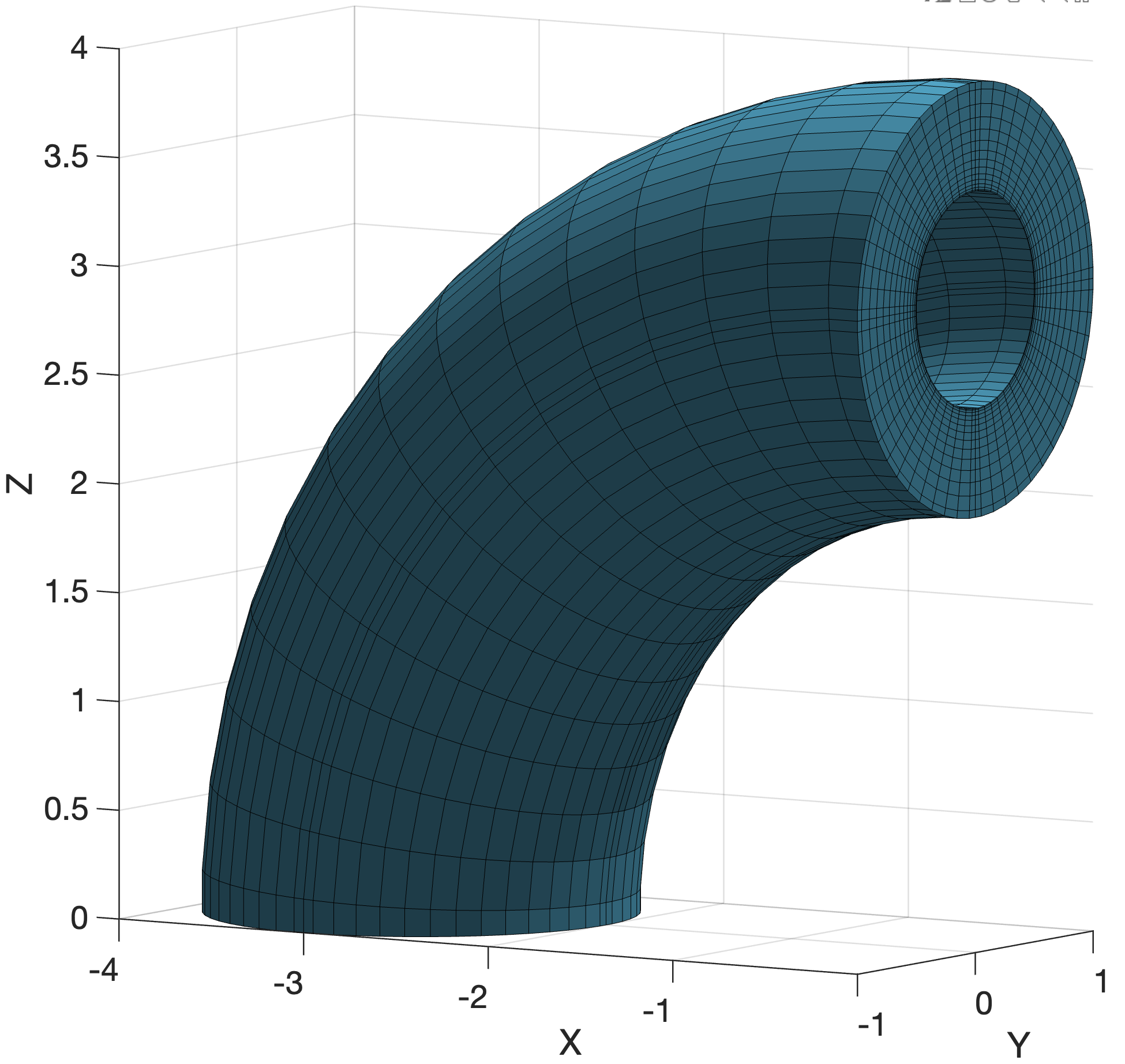}
    \caption{Quarter of torus: $r_\text{in} = 0.5$, $r_\text{out} = 1.0$ and toroidal radius $R = 3$.}
  \end{subfigure}
  \caption{Illustration of control meshes for the closed hemisphere, opened hemisphere, ring, L-shaped, quarter torus, and hyperboloid structures used in the benchmark simulations.}
  \label{Fig.6Geo}
\end{figure}

To further assess the generality and robustness of the TT-IGA framework, we conduct benchmark tests on six distinct three-dimensional geometries: a closed hemisphere, an opened hemisphere, a 3D ring, an L-shaped domain, a hyperboloid, and a quarter torus. The geometric configurations and their control meshes are illustrated in \Figref{Fig.6Geo}. For all cases, the Poisson equation is solved with a distributed source term $f(x, y, z) = \sin(\pi x) \sin(\pi y) \sin(\pi z)$. Dirichlet boundary conditions are imposed as follows: $u = 0$ is applied on the L-shaped domain at $r = r_{\text{in}}$ and $r = r_{\text{out}}$; on the hyperboloid structure at $x = 0$ and $x = H$; and on the open hemisphere at $x = 0$ and $x = 2.0$. 

\begin{figure}
\begin{center}
\minipage{0.45\textwidth}
\includegraphics[width=0.99\textwidth]{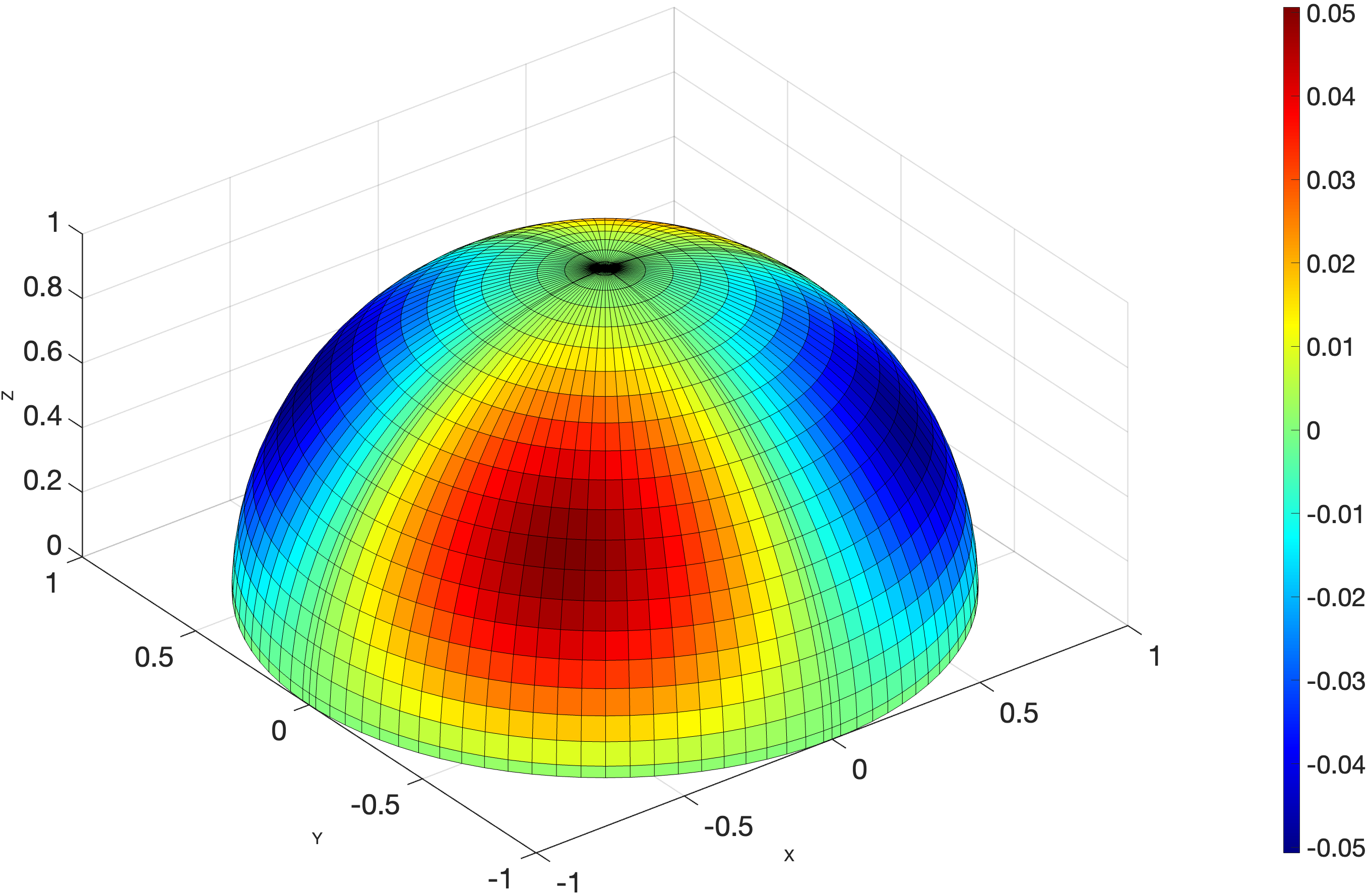}
\caption*{(a) Closed hemisphere}
\endminipage
\hfill
\minipage{0.45\textwidth}
\includegraphics[width=0.99\textwidth]{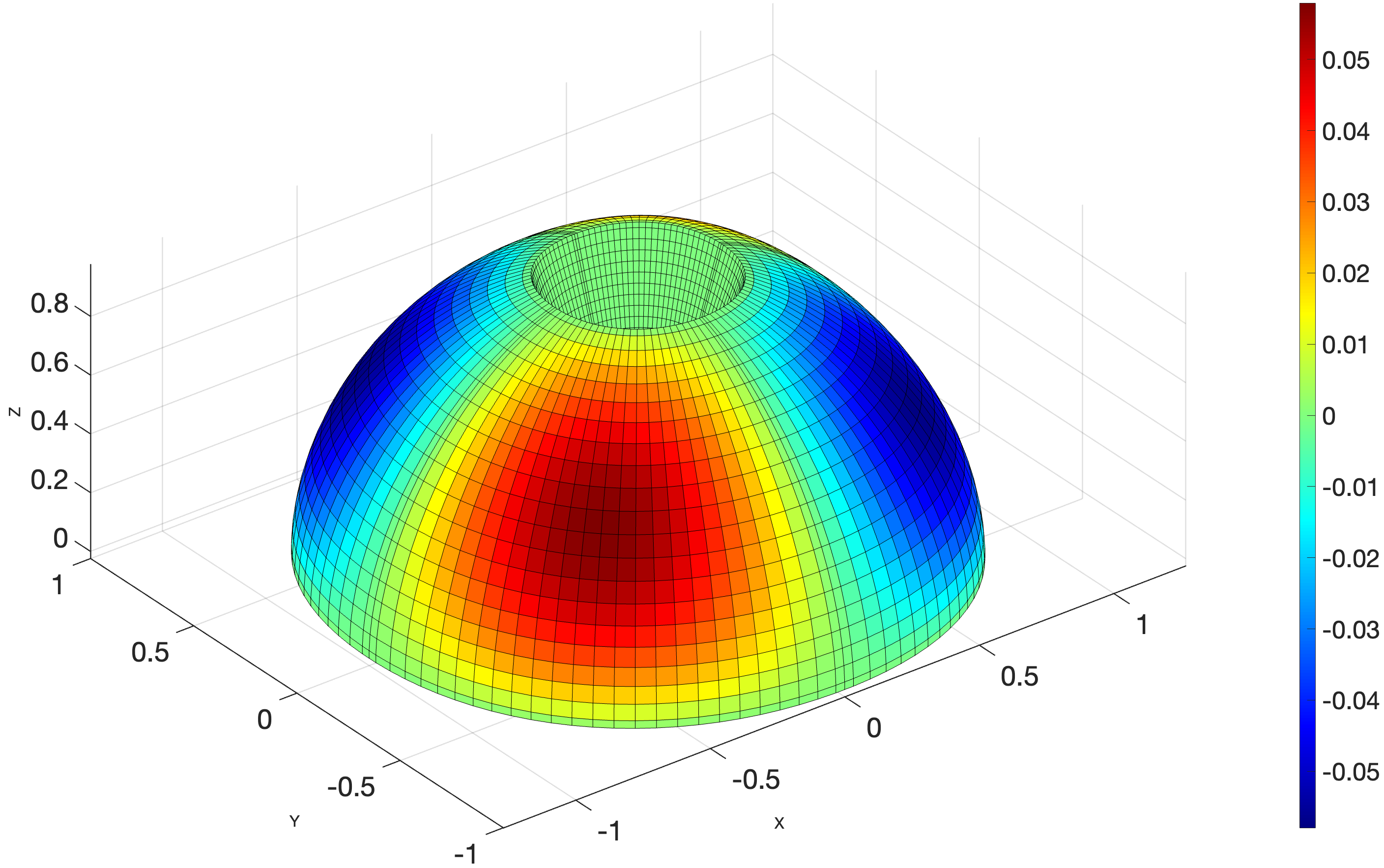}
\caption*{(b) Opened hemisphere}
\endminipage
\vfill
\minipage{0.45\textwidth}
\includegraphics[width=0.99\textwidth]{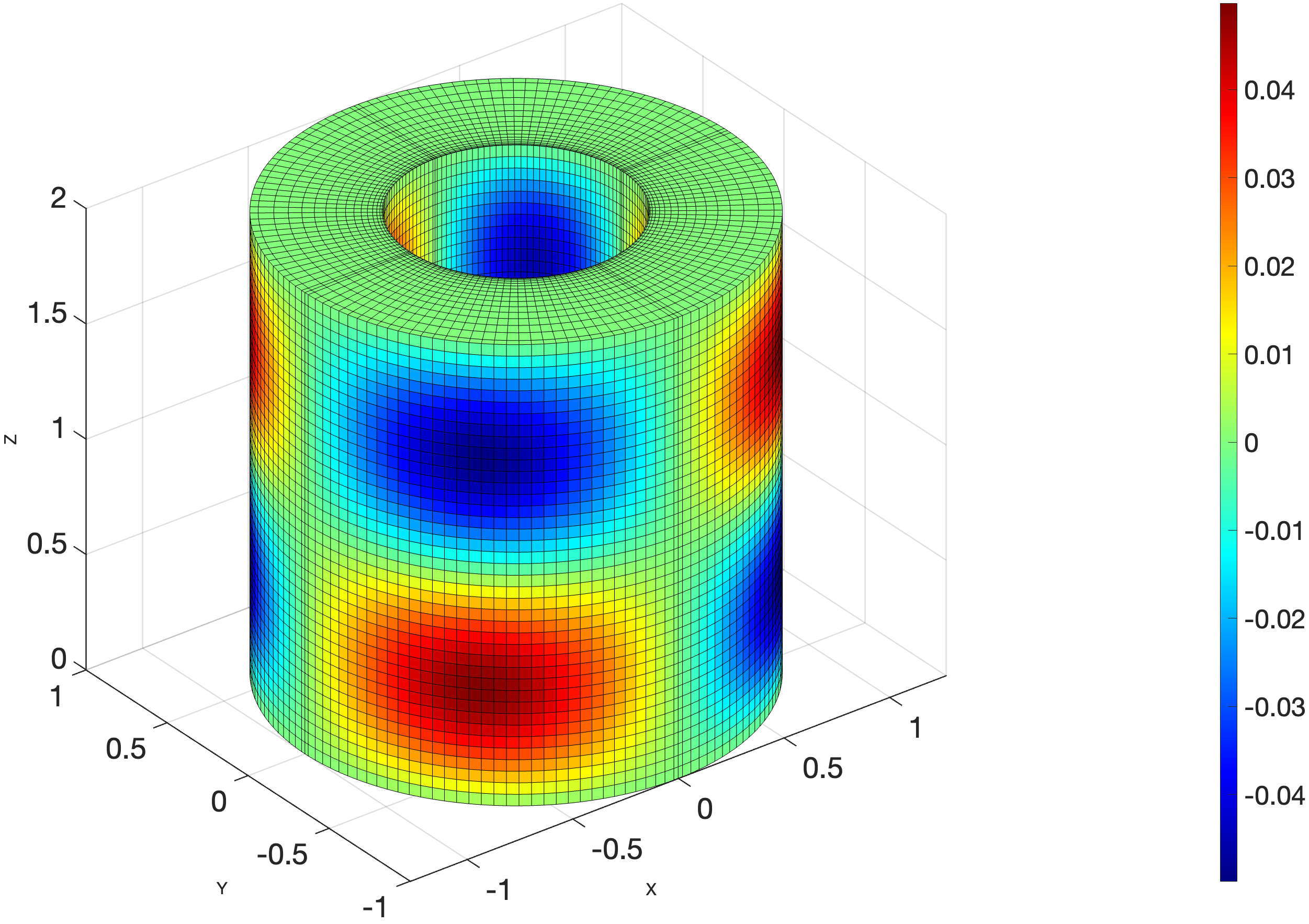}
\caption*{(c) 3D Ring}
\endminipage
\hfill
\minipage{0.45\textwidth}
\includegraphics[width=0.99\textwidth]{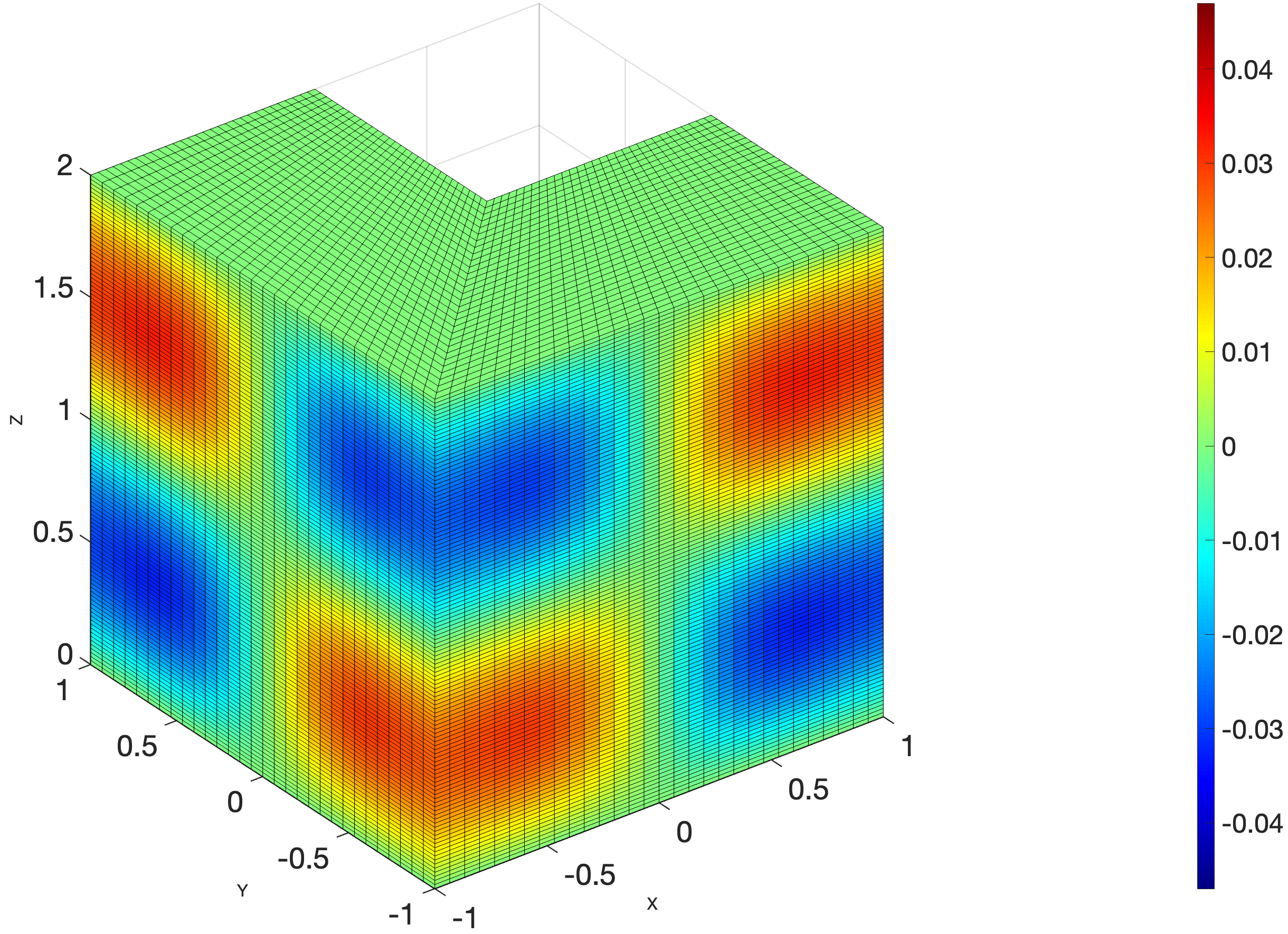}
\caption*{(d) L-shape}
\endminipage
\vfill
\minipage{0.45\textwidth}
\includegraphics[width=0.99\textwidth]{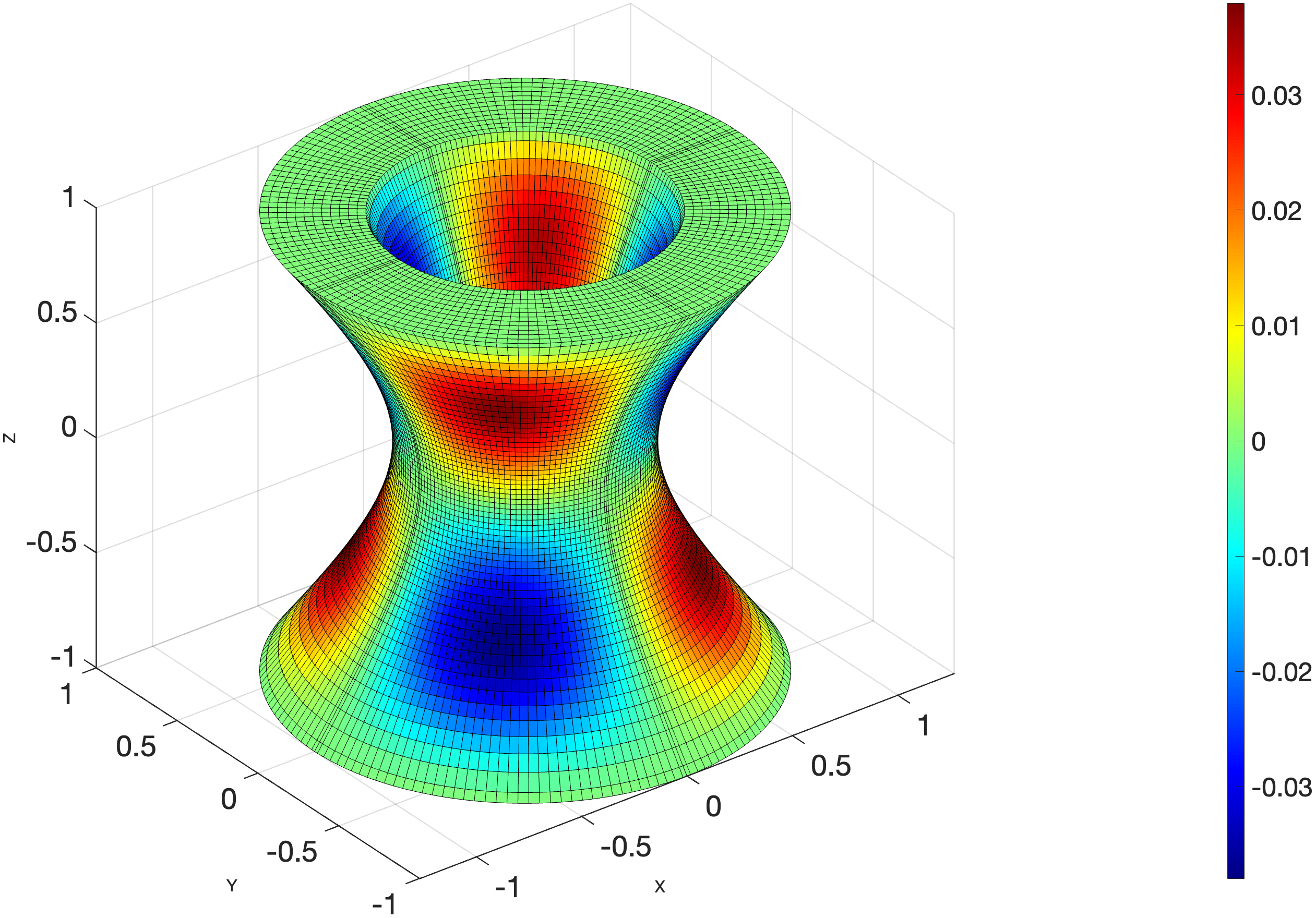}
\caption*{(e) Hyperboloid}
\endminipage
\hfill
\minipage{0.45\textwidth}
\includegraphics[width=0.99\textwidth]{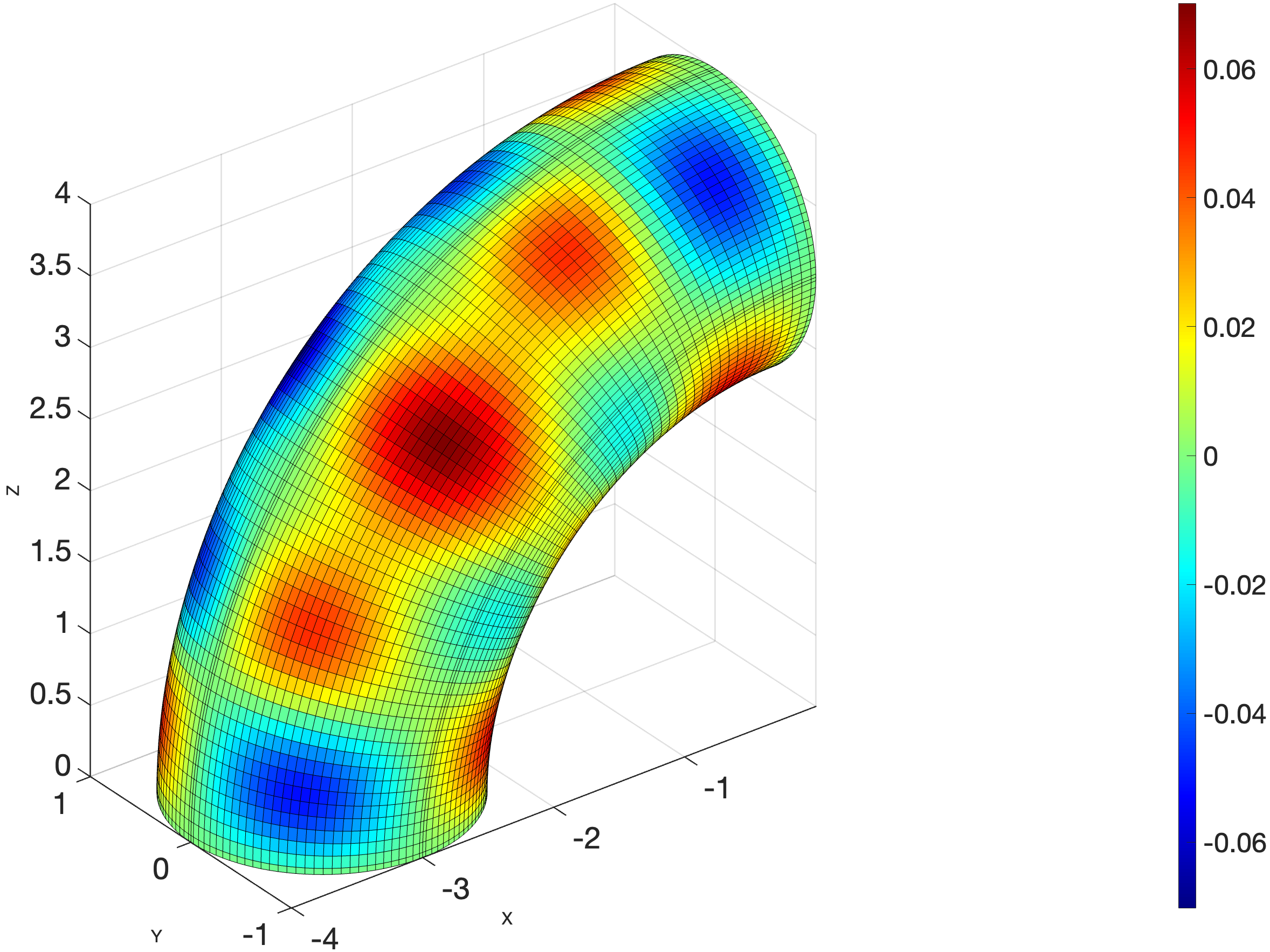}
\caption*{(f) Torus}
\endminipage
\caption{Numerical solutions for six different geometries under the source term $f(x, y, z) = \sin(\pi x) \sin(\pi y) \sin(\pi z)$.}
\label{Fig.6Geo_field}
\end{center}
\end{figure}

\begin{figure}
\begin{center}
\minipage{0.45\textwidth}
\includegraphics[width=0.99\textwidth]{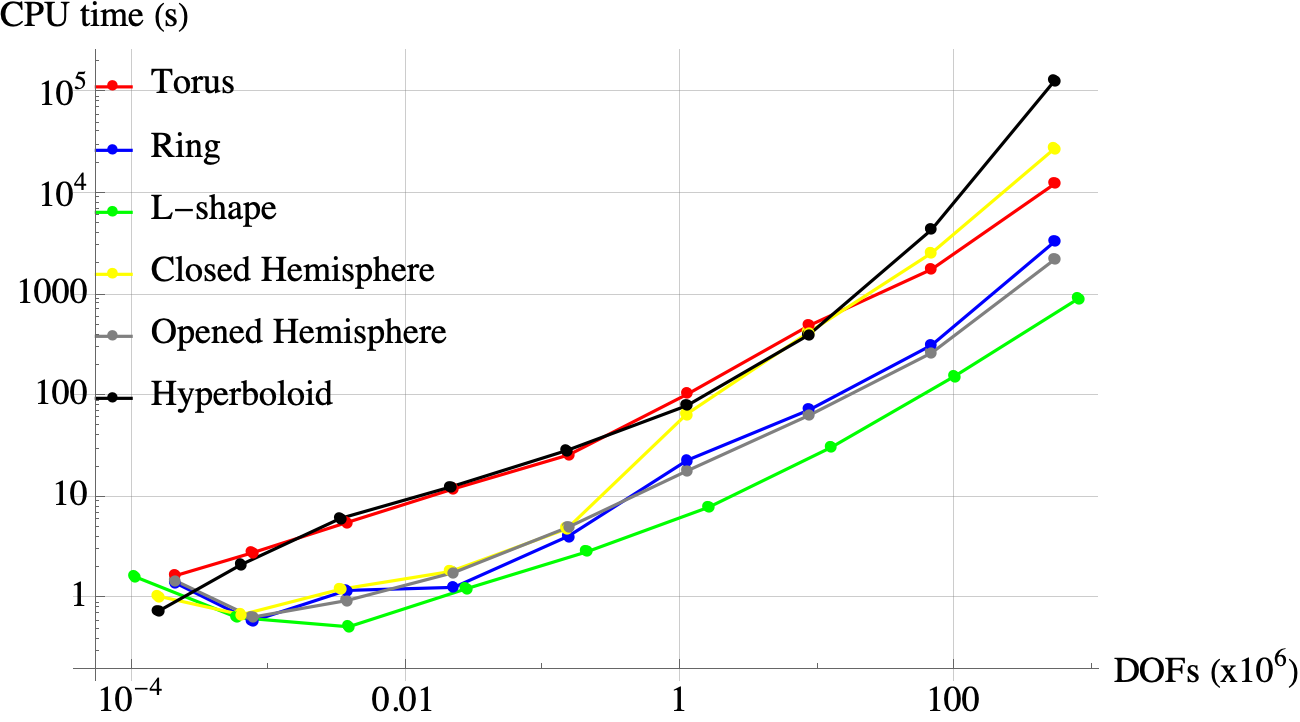}
\caption*{(a) CPU times}
\endminipage
\hfill
\minipage{0.45\textwidth}
\includegraphics[width=0.99\textwidth]{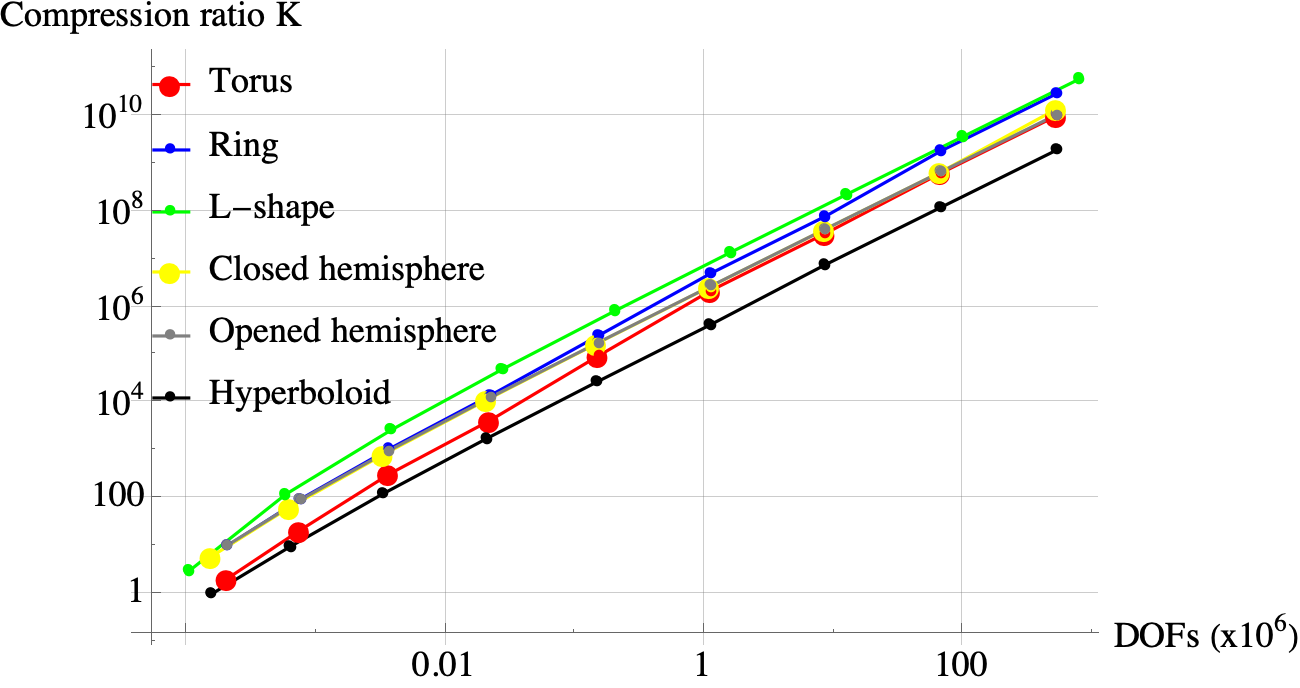}
\caption*{(b) Compression ratio of stiffness matrix K}
\endminipage
\vfill
\minipage{0.45\textwidth}
\includegraphics[width=0.99\textwidth]{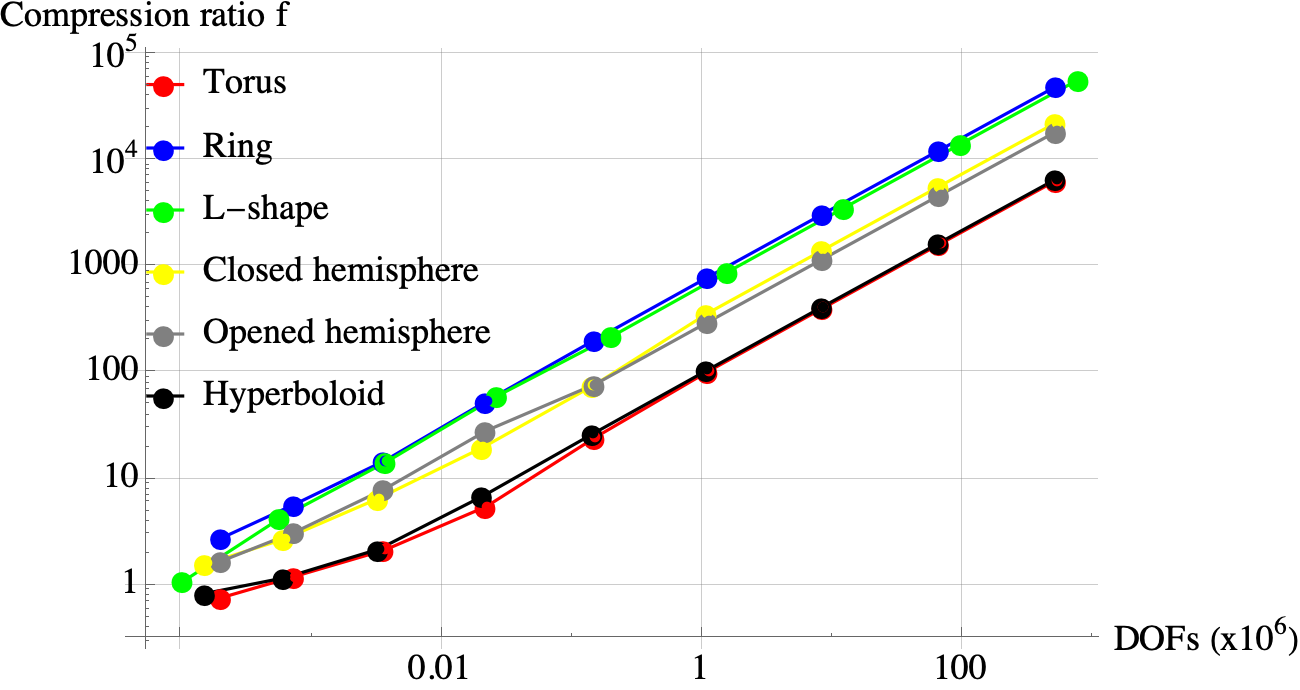}
\caption*{(c) Compression ratio of applied force f}
\endminipage
\hfill
\minipage{0.45\textwidth}
\includegraphics[width=0.99\textwidth]{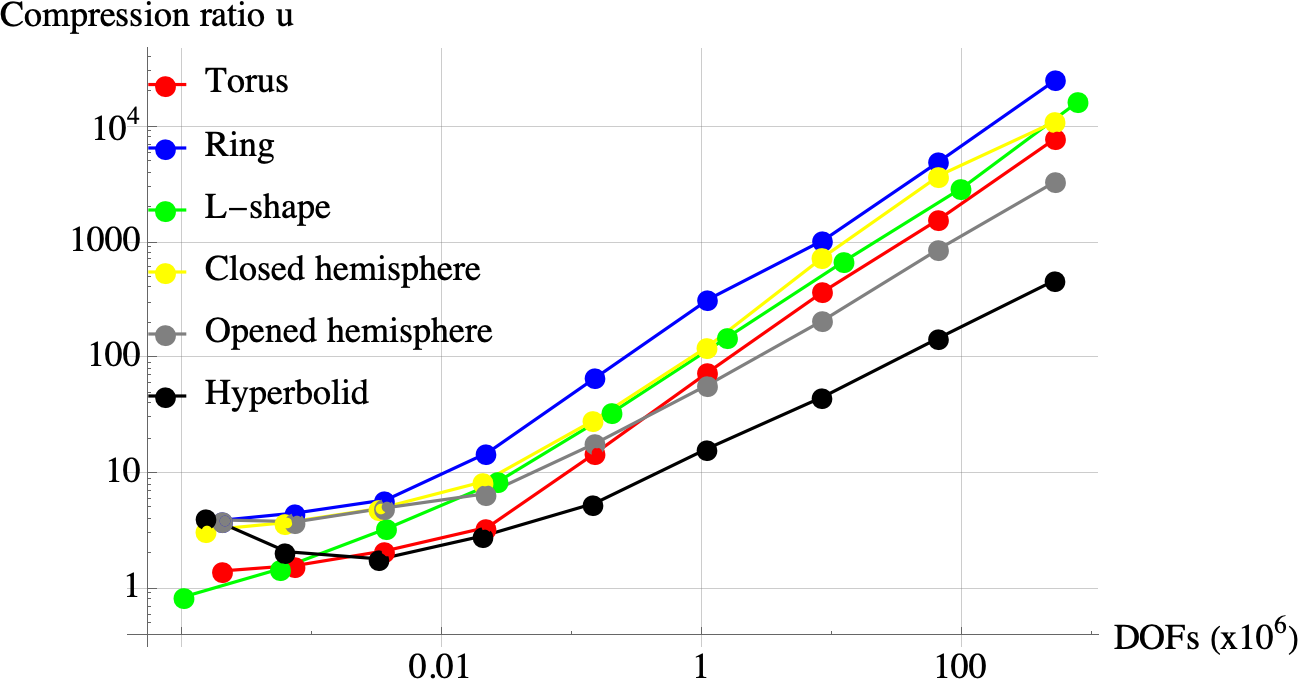}
\caption*{(d) Compression ratio of solution u}
\endminipage
\caption{CPU time and compression ratios of the stiffness matrix $\mathbf{K}$, force vector $\mathbf{f}$, and solution $\mathbf{u}$ for the six tested geometries.}
\label{Fig.6Geo_graph}
\end{center}
\end{figure}

\Figref{Fig.6Geo_field} displays representative numerical solutions for each geometry, demonstrating the flexibility of TT-IGA in handling several complex domains. For each test, the number of degrees of freedom is increased up to approximately 0.5 billion, showcasing the scalability of the approach. Performance metrics, including CPU time and compression ratios for the stiffness matrix $\mathbf{K}$, force vector $\mathbf{f}$, and solution vector $\mathbf{u}$, are summarized in  \Figref{Fig.6Geo_graph}. The results show consistently excellent compression across all tested geometries: as the problem size increases, both memory usage and computational time grow much more slowly than in traditional full-grid methods. Remarkably, the TT-IGA method maintains its efficiency and accuracy for all tested shapes, including the most challenging cases such as the hyperboloid and quarter torus.


\section{Conclusion}
\label{sec::conclude}

We have introduced TT-IGA, a tensor-train-based isogeometric analysis framework for solving the 3D Poisson equation on complex geometries. The method combines the geometric flexibility and higher-order accuracy of IGA with the compression and scalability of tensor network solvers. By exploiting the tensor-product structure of spline basis functions and compressing geometry-dependent quantities via TT-cross interpolation, we assemble the stiffness matrix and load vector entirely in TT format without sacrificing geometric fidelity.

Our numerical experiments demonstrate that TT-IGA achieves convergence rates consistent with standard IGA while offering significant memory savings and computational efficiency. The solver remains robust across a variety of geometries—including L-shapes, rings, hemispheres, toroidal and hyperboloidal domains—with degrees of freedom exceeding 500 million. Compression ratios improve with mesh refinement, confirming that TT-IGA is especially effective at large scale.

Unlike earlier TT-IGA formulations that rely on separable or parameterized domains, our framework directly handles exact NURBS-based geometries without approximation. This opens the door to applying TT solvers in realistic, high-fidelity simulations beyond simple Cartesian grids.

In summary, TT-IGA offers a scalable and geometry-aware strategy for solving elliptic PDEs, providing both high accuracy and low computational cost. The proposed framework paves the way for tensor network-based solvers in advanced engineering simulations involving complex 3D domains.
\section*{Acknowledgments}
The authors gratefully acknowledge the support of the Laboratory Directed Research and Development (LDRD) program of Los Alamos National Laboratory under project number 20230067DR and 20250893ER.
Los Alamos National Laboratory is operated by Triad National Security, LLC, for the National Nuclear Security Administration of U.S. Department of Energy (Contract No.\ 89233218CNA000001).

\bibliographystyle{IEEEtran}
\bibliography{References}
\end{document}